\documentclass[12pt]{amsart}

\newtheorem{theorem}{Theorem}[section]
\newtheorem{proposition}[theorem]{Proposition}
\newtheorem{lemma}[theorem]{Lemma}
\newtheorem{corollary}[theorem]{Corollary}
\newtheorem{remark}[theorem]{Remark}
\newtheorem{definition}[theorem]{Definition}

\numberwithin{equation}{section}

\begin{document}
\baselineskip=15pt

\title[Fundamental group scheme and ample hypersurface]{Comparison
of fundamental group schemes of a projective variety and
an ample hypersurface}

\author[I. Biswas]{Indranil Biswas}

\address{School of Mathematics, Tata Institute of Fundamental
Research, Homi Bhabha Road, Bombay 400005, India}

\email{indranil@math.tifr.res.in}

\author[Y. I. Holla]{Yogish I. Holla}

\address{School of Mathematics, Tata Institute of Fundamental
Research, Homi Bhabha Road, Bombay 400005, India}

\email{yogi@math.tifr.res.in}

\date{}

\begin{abstract}

Let $X$ be a smooth projective variety defined over an algebraically
closed field, and let $L$ be an ample line bundle over $X$. We prove that
for any smooth hypersurface $D$ on $X$ in the complete linear system
$\vert L^{\otimes d}\vert$, the inclusion map $D\hookrightarrow X$
induces an isomorphism of fundamental group schemes, provided
$d$ is sufficiently large and $\dim X \, \geq\, 3$. If
$\dim X \, =\, 2$, and $d$ is sufficiently large, then the induced
homomorphism of fundamental group schemes remains surjective. We give
an example to show that the homomorphism of
fundamental group schemes induced by the inclusion map of a
reduced ample curve in a smooth projective surface
is not surjective in general.

\end{abstract}

\maketitle

\section{Introduction}

Let $X$ be a projective variety defined over an
algebraically closed field $k$. Let $D$ be a
reduced ample hypersurface in $X$. If $X$ is normal, and
$\dim X \,\geq\, 2$, then it is known that the
induced homomorphism of \'etale fundamental groups
$$
\pi_1(D,x_0)\, \longrightarrow\, \pi_1(X,x_0)
$$
is surjective; here $x_0$ is any $k$--rational point
of $D$. If $X$ is smooth with $\dim X\, \geq\, 3$, then from
Grothendieck's Lefschetz theory it follows that the above
homomorphism of \'etale fundamental groups is an
isomorphism (see \cite[Expos\'e X]{Gr}).

When the characteristic of $k$ is positive, Nori
constructed an invariant of $X$ which is
finer than the \'etale fundamental group $\pi_1(X,x_0)$
\cite{No1}. He calls this invariant
the fundamental group scheme of $X$.

To describe the fundamental group scheme, we first recall that
the \'etale fundamental group $\pi_1(X,x_0)$ coincides
with group scheme associated to the Tannakian category
defined by all vector bundles $E$ over $X$ with the following
property: there is a finite \'etale Galois cover $Y$ of $X$
such that the pull back of $E$ to $Y$ is trivializable. The
fiber functor for this Tannakian category sends $E$
to its fiber $E_{x_0}$. The fundamental group scheme
is the group scheme associated to the Tannakian category
defined by all vector bundles $E$ over $X$ with the following
property: there is a principal bundle $Y\, \longrightarrow\, X$
over $X$
with a finite group scheme as the structure group scheme
such that the pull back of $E$ to $Y$ is trivializable.
Note that if the characteristic of the base field $k$ is zero,
then the fundamental group scheme coincides with the
\'etale fundamental group.

The fundamental group scheme of $X$ with
$x_0$ as the base point is denoted by
$\pi(X,x_0)$. Comparing the above definitions
of $\pi_1(X,x_0)$ and $\pi(X,x_0)$ it follows
that there is a canonical surjective (faithfully
flat) homomorphism from
$\pi(X,x_0)$ to $\pi_1(X,x_0)$. The kernel of this
surjective homomorphism $\pi(X,x_0)\,\longrightarrow\,
\pi_1(X,x_0)$ is a local group scheme.

Although in the case of positive characteristics the group
scheme $\pi(X,x_0)$ can be larger than $\pi_1(X,x_0)$, 
the fundamental group scheme continues to share
some of the basic properties of the \'etale fundamental group.
Our aim here is to investigate the relationship 
between $\pi(D,x_0)$ and $\pi(X,x_0)$, where $D\, \subset\, X$
is an ample hypersurface.

We first give an example of a pair $(X\, ,D)$, where $X$ is a
smooth projective surface and $D$ a reduced ample
curve in $X$, such that the natural homomorphism
$$
\pi(D,x_0)\, \longrightarrow\, \pi(X,x_0)
$$
is not surjective. This example is based on the fact
that the Kodaira vanishing theorem may fail in positive
characteristics. More precisely, given an ample
hypersurface $D\, \hookrightarrow\, X$ with $\dim X\,
\geq\, 2$, the induced homomorphism
$$
\pi(D,x_0)\, \longrightarrow\, \pi(X,x_0)
$$
fails to be surjective whenever
$H^1(X,\, {\mathcal O}_X(-D))\, \not=\, 0$.

On the other hand, the following theorem shows that
the above homomorphism between fundamental group
schemes behaves like the homomorphism between
\'etale fundamental groups if $D$ is sufficiently
positive.

\begin{theorem}\label{thm0}
Let $X$ be a smooth projective variety defined over an
algebraically closed field and $L$ an ample line bundle over $X$.

\begin{enumerate}
\item{} Assume that $\dim X \, \geq\, 2$. There is an integer
$d_1(X,L)$ such that for any smooth divisor
$D\, \in\, \vert L^{\otimes d}\vert$, where $d\, >\, d_1(X,L)$,
the natural homomorphism between fundamental group schemes
\begin{equation}\label{1}
\rho\, :\, \pi(D,x_0) \,\longrightarrow\, \pi(X,x_0)
\end{equation}
is surjective (faithfully
flat), where $x_0$ is any $k$--rational point of $D$.

\item{} Assume that $\dim X \, \geq\, 3$. There is an integer
$d_2(X,L)$ such that for any smooth divisor
$D\, \in\, \vert L^{\otimes d}\vert$, where $d\, >\, d_2(X,L)$,
the homomorphism $\rho$ in Eqn. \eqref{1} is a closed immersion.
(So combining with statement (1) it follows that $\rho$ is an
isomorphism if  $d\, >\, d_1(X,L), d_2(X,L)$.)
\end{enumerate}
\end{theorem}

An effective estimate of the constant $d_1(X,L)$ in
Theorem \ref{thm0}(1) is given in Theorem
\ref{thm1}. An effective estimate of the constant $d_2(X,L)$ in
the second part of the theorem
is described in Remark \ref{inj.eff.}.

The first part of the above theorem is proved in Theorem
\ref{thm1}, and the second part is proved in Theorem \ref{thm2}.

In Section \ref{sec.ex}, the earlier mentioned example is
constructed. In Section \ref{sec.s.v.r.}, some results needed
in the proof of Theorem \ref{thm2} are established.

\section{An example}\label{sec.ex}

Let $k$ be an algebraically closed field.
Let $F_k$ be the Frobenius homomorphism of the field $k$.
Let $X$ be a variety defined over $k$.
We then have the geometric Frobenius morphism
$$
F_X\, :\, X\, \longrightarrow\, F^*_k X\, .
$$
Note that since the field $k$
is perfect, the variety $F^*_k X$ is isomorphic to
$X$. For notational
convenience, we will denote the
iterated pull back $(F^n_k)^* X$
by $X$ itself, where $n$ is any positive integer.
With this notation, the $n$--fold iteration of
the geometric Frobenius morphism
$$
F^n_X \, :\, X\, \longrightarrow\, (F^n_k)^* X
$$
will be considered as a morphism from $X$ to $X$.
We will
use this convention throughout in the manuscript.

\begin{definition}\label{def1}
{\rm A vector bundle $E$ over $X$ will be called
an} $F$--\textit{trivial vector bundle}
{\rm if there is a nonnegative integer $m$ such that
the vector bundle
$(F^m_X)^*E$ is isomorphic to a trivial vector bundle over $X$.}
\end{definition}

Let $X$ be a projective variety defined over $k$.
For any $k$--rational
point $x_0\, \in\, X$, the \'etale fundamental group of $X$
with $x_0$ as the base point will be denoted by $\pi_1(X,x_0)$.

Let $D\,\subset\, X$ be a reduced ample
hypersurface in $X$. If we have $\dim X\, \geq\, 2$
with $X$ normal, then from
Grothendieck's Lefschetz theory, \cite[Expos\'e X]{Gr}, it follows
that for
any $k$--rational point $x_0\, \in\, D$, the natural homomorphism
\begin{equation}\label{r1}
\rho_1\, :\, \pi_1(D,x_0)\, \longrightarrow\, \pi_1(X,x_0)
\end{equation}
is surjective. To prove that $\rho_1$ is surjective, take any
connected \'etale Galois cover
$$
\phi\, :\, \widetilde{X}\, \longrightarrow\, X\, .
$$
Since $X$ is normal, we conclude that $\widetilde{X}$ is irreducible.
As the divisor is reduced and ample, the inverse image
$\phi^{-1}(D)$ is also a reduced and ample hypersurface in
$\widetilde{X}$. Therefore, we know that $\phi^{-1}(D)$ is 
connected (see \cite[page 79, Corollary 6.2]{Ha}
and \cite[page 64, Proposition 2.1]{Ha}). This immediately implies
that the homomorphism $\rho$ in Eqn. \eqref{r1} is surjective.

For any $k$--rational point $x_0$ of $D$, the fundamental
group scheme of $X$ (respectively, $D$) with $x_0$ as the base
point will be
denoted by $\pi(X,x_0)$ (respectively, $\pi(D,x_0)$); the
fundamental group scheme was introduced in \cite{No1}, \cite{No2}.

Let
\begin{equation}\label{rho}
\rho\, :\, \pi(D,x_0)\, \longrightarrow\, \pi(X,x_0)
\end{equation}
be the homomorphism of fundamental group schemes induced by the
inclusion map $D\, \hookrightarrow\, X$. Our aim
in this section is to give
an example of a pair $(X\, ,D)$, where $D$ is
an ample reduced curve on a smooth projective
surface $X$, for which the homomorphism
$\rho$ in Eqn. \eqref{rho} is not surjective.
This is in contrast with the situation for
\'etale fundamental groups.

Assume that the characteristic of the field $k$ is positive.
In \cite{Ra}, Raynaud constructed a pair $(S\, , C)$,
where $S$ is a smooth projective surface and
$C$ a reduced ample hypersurface on $S$, such that
\begin{equation}\label{nz}
H^1(S,\, {\mathcal O}_{S}(-C)) \, \not=\, 0\, .
\end{equation}
It may be noted that the Kodaira vanishing theorem says that
$H^1(X', \,{\mathcal O}_{X'}(-D'))\,=\, 0$ for any ample divisor
$D'$ on a smooth projective surface $X'$ defined over an
algebraically closed field of characteristic zero.

For our example, set $X \, =\, S$, and take $D$
to be the divisor $C$ on $S$ in the above mentioned
example $(S\, , C)$ of \cite{Ra}.

We consider the natural short exact sequence of sheaves
$$
0 \, \longrightarrow\, {\mathcal O}_X(-D)
\, \longrightarrow\,{\mathcal O}_X\, \longrightarrow\,
{\mathcal O}_D\, \longrightarrow\, 0
$$
on $X$. Let
\begin{equation}\label{esc}
H^0(X,\, {\mathcal O}_X) \, =\, H^0(D,\, {\mathcal O}_D)
\,\longrightarrow\, H^1(X,\, {\mathcal O}_X(-D)) \,
\stackrel{\gamma}{\longrightarrow}\, H^1(X,\, {\mathcal O}_X)\,
\stackrel{\delta}{\longrightarrow}\, H^1(D,\, {\mathcal O}_D)
\end{equation}
be a segment of the long exact sequence of
cohomologies corresponding
to this exact sequence of sheaves. Take any nonzero element
\begin{equation}\label{defc}
c\, \in\, H^1(X,\, {\mathcal O}_X(-D))\setminus \{0\}
\end{equation}
which exists by Eqn. \eqref{nz}. Therefore, the
cohomology class $c$ gives a
non--split exact sequence of vector bundles
\begin{equation}\label{fex}
0\, \longrightarrow\, {\mathcal O}_X \, \longrightarrow\, W
\, \stackrel{\psi}{\longrightarrow}\,{\mathcal O}_X(D) \,
\longrightarrow\, 0
\end{equation}
over $X$.

Set
\begin{equation}\label{alpha}
\alpha\, :=\, \gamma(c)\, \in\, H^1(X,\, {\mathcal O}_X)
\setminus\{0\}\, ,
\end{equation}
where $\gamma$ is the homomorphism in Eqn. \eqref{esc} and $c$
is the element in Eqn. \eqref{defc}. Let
\begin{equation}\label{V}
0\, \longrightarrow\, {\mathcal O}_X \, \longrightarrow\, V
\, \longrightarrow\,{\mathcal O}_X \, \longrightarrow\, 0
\end{equation}
be the nontrivial extension corresponding to the cohomology
class $\alpha$ defined in Eqn. \eqref{alpha}. So we have
$V\, =\, \psi^{-1}({\mathcal O}_X)$, where $\psi$
is the projection in Eqn. \eqref{fex}.

Since
$$
\delta (\alpha)\, =\, \delta (\gamma(c)) \, =\, 0\, ,
$$
where $\delta$ is the homomorphism in Eqn. \eqref{esc},
we conclude that the
restriction of the vector bundle $V$ to $D$ splits as
\begin{equation}\label{Vr}
V\vert_D \, =\, {\mathcal O}_D\oplus {\mathcal O}_D\, .
\end{equation}

Let $F_X\, :\, X\, \longrightarrow\, X$ be the
Frobenius morphism. For any positive integer $n$, let
$$
F^n_X \, :=\, \overbrace{F_X\circ\cdots\circ F_X}^{n\mbox{-}\rm{times}}
\, :\, X\, \longrightarrow\, X
$$
be the $n$--fold iteration of the self--morphism $F_X$. By $F^0_X$
we will mean the identity morphism of $X$.

\begin{lemma}\label{lem1}
There is a positive integer $n$ such that the
vector bundle $(F^n_X)^* V$ over $X$ is isomorphic
to ${\mathcal O}_X\oplus {\mathcal O}_X$, where $V$ is
the vector bundle constructed in Eqn. \eqref{V}.
\end{lemma}

\begin{proof}
For any morphism $f\, :\, Y\, \longrightarrow\, X$,
consider the short exact sequence of vector bundles over $Y$
$$
0\, \longrightarrow\, f^*{\mathcal O}_X \, =\, {\mathcal O}_Y\,
\longrightarrow\, f^*V \, \longrightarrow\,
f^*{\mathcal O}_X \, =\, {\mathcal O}_Y\, \longrightarrow\, 0
$$
obtained by pulling back the exact sequence
Eqn. \eqref{V}. The cohomology class
in $H^1(Y,\, {\mathcal O}_Y)$ corresponding
to it coincides with
$$
f^*\alpha \, \in\, H^1(Y,\, f^*{\mathcal O}_X) \, =\,
H^1(Y,\, {\mathcal O}_Y)\, ,
$$
where $\alpha$ is defined in Eqn. \eqref{alpha}. Let
\begin{equation}\label{defh}
h\, :\, H^1(Y,\, {\mathcal O}_Y(-f^{-1}(D)))\, \longrightarrow\,
H^1(Y,\, {\mathcal O}_Y)
\end{equation}
by the homomorphism obtained from the exact sequence
$$
0 \, \longrightarrow\, {\mathcal O}_Y(-f^{-1}(D))
\, \longrightarrow\,{\mathcal O}_Y\, \longrightarrow\,
{\mathcal O}_{f^{-1}(D)}\, \longrightarrow\, 0\, .
$$
{}From the definition of $\alpha$ it follows that
\begin{equation}\label{1in}
f^*\alpha \, \in\, h(H^1(Y,\, {\mathcal O}_Y(-f^{-1}(D))))\, ,
\end{equation}
where $h$ is the homomorphism in Eqn. \eqref{defh}. Indeed, the
identity $f^*\alpha\, =\, h(f^*c)$ holds, 
where $c$ is the cohomology class in Eqn. \eqref{defc}.

We set $f \, :=\, F^m_X\, :\, X\longrightarrow\, X$. Then
\begin{equation}\label{pid}
f^{-1}(D) \, =\, p^m D\, ,
\end{equation}
where $p$ is the characteristic of the
field $k$ (which is positive by our assumption).

Let $K_X$ denote the canonical line bundle of $X$.
Since $D$ is an ample divisor on $X$, and
$\dim X \, =\, 2$, we have
$$
H^1(X,\, {\mathcal O}_X(-jD)) \,=\,
H^1(X,\, K_X\otimes {\mathcal O}_X(jD))^*\, =\, 0
$$
for all $j$ sufficiently large
\cite[page 111, Proposition 6.2.1]{EGA3}, where the
first isomorphism is the Serre duality.
Therefore, using Eqn. \eqref{pid} we conclude that
$$
H^1(X,\, {\mathcal O}_X(-f^{-1}(D)))\, =\, 0
$$
for all $m$ sufficiently large, where
$f \, =\, F^m_X$. Consequently, from Eqn. \eqref{1in}
it follows that $(F^m_X)^*\alpha \, =\, 0$ for all
$m$ sufficiently large. Hence the exact sequence
$$
0\, \longrightarrow\,{\mathcal O}_X\,\longrightarrow\, (F^m_X)^*V
\, \longrightarrow\, {\mathcal O}_X\, \longrightarrow\, 0
$$
obtained by pulling back the exact sequence in Eqn. \eqref{V}
using $F^m_X$
splits for all $m$ sufficiently large. This completes the
proof of the lemma.
\end{proof}

We will now prove a proposition from which it will follow
that the vector bundle $V$ in Eqn. \eqref{V}
is essentially finite; see \cite[page 38]{No1} for the
definition of essentially finite vector bundles.

\begin{proposition}\label{pr.e.f.}
Let $M$ be a projective variety defined over an algebraically
closed field $k$. Let $E$ be a vector bundle over $M$ with the
following property. There is an \'etale Galois covering
$$
\beta\, :\, \widetilde{M}\, \longrightarrow\, M
$$
such that the vector bundle $\beta^*E$ over
$\widetilde{M}$ is $F$--trivial (see Definition \ref{def1}).
Then $E$ is an essentially finite vector bundles.
\end{proposition}

\begin{proof}
Let $r$ denote the rank of $E$. Let $n$ be a positive integer
such that the vector bundle $(F^n_{\widetilde{M}})^*\beta^*E$ is
trivializable, where
$$
F_{\widetilde{M}}\, :\,
\widetilde{M}\, \longrightarrow\, \widetilde{M}
$$
is the Frobenius morphism.

Consider the exact sequence of group schemes
$$
e\, \longrightarrow\, \text{kernel}(F^n_{\text{GL}(r,k)})
\, \longrightarrow\, \text{GL}(r,k) \,
\stackrel{F^n_{\text{GL}(r,k)}}{\longrightarrow}\,
\text{GL}(r,k) \, \longrightarrow\, e\, ,
$$
where $F_{\text{GL}(r,k)}\, :\, \text{GL}(r,k) \, \longrightarrow\,
\text{GL}(r,k)$ is the Frobenius morphism. Let
\begin{equation}\label{13d2}
H^1(\widetilde{M}, \, \text{kernel}(F^n_{\text{GL}(r,k)}))
\, \longrightarrow\, H^1(\widetilde{M}, \,\text{GL}(r,k))
\, \stackrel{\phi}{\longrightarrow}\,
H^1(\widetilde{M}, \,\text{GL}(r,k))
\end{equation}
be the exact sequence of pointed sets
obtained from this short exact sequence.

For any
principal $G$--bundle $E_G$ over $\widetilde{M}$,
where $G$ is an algebraic group defined over $k$,
the pull back
$F^*_{\widetilde{M}} E_G$ is identified with the principal
$G$--bundle over $\widetilde{M}$
obtained by extending the structure group of $E_G$ using the
Frobenius morphism $F_G\, :\, G\, \longrightarrow\, G$
\cite[287, Remark 3.22]{RR}. Using this, together with the
fact that $(F^n_{\widetilde{M}})^*\beta^*E$ trivializable,
we conclude that the element $\alpha\, \in\,
H^1(\widetilde{M}, \,\text{GL}(r,k))$ corresponding to
$(F^n_{\widetilde{M}})^*\beta^*E$ has the property that
$\phi(\alpha)$ is the base point in
$H^1(\widetilde{M}, \,\text{GL}(r,k))$ (the point corresponding to
the trivial $\text{GL}(r,k)$--bundle), where $\phi$
is the map in Eqn. \eqref{13d2}.

This implies that the vector bundle $\beta^*E$
is associated to a principal bundle
over $\widetilde{M}$ with the finite group scheme
$\text{kernel}(F^n_{\text{GL}(r,k)})$ as the
structure group scheme. Since $\beta$ is an \'etale
Galois covering, from this it follows that
the vector bundle $E$
is associated to a principal bundle over $M$
with the finite group scheme
as the structure group scheme. Consequently,
the vector bundle $E$ is essentially finite
\cite[page 38, Proposition 3.8]{No1}.
This completes the proof of the proposition.
\end{proof}

Using Proposition \ref{pr.e.f.}, from
Lemma \ref{lem1} it follows that the vector
bundle $V$ in Eqn. \eqref{V}
is essentially finite.
The vector bundle $V$, being essentially finite, corresponds
to a representation of the fundamental group scheme
$\pi(X,x_0)$, where $x_0$ is any $k$--rational point of $X$.
The fundamental group scheme is defined in \cite[page 40]{No1}.

Consider the two essentially
finite vector bundles over $X$, namely $V$ and
${\mathcal O}_X\oplus {\mathcal O}_X$. Since $V$
is a nontrivial extension of ${\mathcal O}_X$ by ${\mathcal O}_X$,
these two vector bundles are not isomorphic. On the other
hand, Eqn. \eqref{Vr} says that
their restrictions to $D$ are isomorphic.

Consequently, we have two non-isomorphic
representations of $\pi(X,x_0)$,
namely $V$ and ${\mathcal O}_X\oplus {\mathcal O}_X$, with the
following property: when these two are considered as representation
of $\pi(D,x_0)$ using $\rho$ in Eqn. \eqref{rho}, then they become
isomorphic representations of $\pi(D,x_0)$.

{}From this it follows immediately that the
homomorphism $\rho$ (defined in Eqn. \eqref{rho}) for this pair
$(X\, ,D)$ is not surjective.

\begin{remark}\label{rem1}
{\rm Let $X$ be a projective variety of dimension at least two.
Let $D\, \subset\, X$ be an ample hypersurface such that
$H^1(X,\, {\mathcal O}_X(-D))\, \not=\, 0$. Then the
above arguments show that the induced homomorphism
$\pi(D,x_0)\, \longrightarrow\, \pi(X,x_0)$ is not surjective.
See \cite[page 120, Proposition 2.14]{Ek} for further examples
of such pairs $(X\, , D)$.}
\end{remark}

\begin{remark}
{\rm Consider the vector bundle $W$ in Eqn. \eqref{fex}.
We note that $c_1(W) \, =\, D$ and $c_2(W) \, =\, 0$. Since
$D$ is an ample divisor, we know that $D^2\, > \, 0$. Therefore,
$c_1(W)^2 \, >\, 4c_2(W) \,=\, 0$. In other words,
the vector bundle $W$ violates the Bogomolov inequality
condition. However the vector bundle $W$
is semistable \cite[pages 251--252]{Mu}. More precisely,
in \cite[pages 251--252]{Mu} Mumford shows that if $W$ is not
semistable, then the exact sequence in Eqn. \eqref{fex} splits.
Hence $W$ provides a counter-example to the Bogomolov inequality
for positive characteristic.}
\end{remark}

\section{Fundamental group scheme and hyperplane section}\label{se2}

Let $k$ be an algebraically closed field of characteristic $p$,
with $p\, >\, 0$.

Let $Y$ be a smooth projective variety defined over $k$. Fix an
ample line bundle $L$ over $Y$ to define degree of coherent
sheaves on $Y$. Let
$$
F_Y\, :\, Y\, \longrightarrow\, Y
$$
be the Frobenius morphism of $Y$.

Take any vector bundle $E$ over $Y$. We will
briefly recall the definition of $L_{\rm max}(E)$
(see \cite[page 257]{La}). For any positive integer $m$, let
$$
0\, =\, E_{m,0} \, \subset\, E_{m,1} \, \subset\, \cdots \, \subset\,
E_{m, a(m)-1} \, \subset\, E_{m, a(m)}\, =\, (F^m_Y)^*E
$$
be the Harder--Narasimhan filtration of the vector bundle $(F^m_Y)^*E$
for the polarization $L$ on $Y$.
It is known that for any sufficiently large $m$, we have
\begin{equation}\label{isop}
E_{m+n,1}\, =\, (F^n_Y)^* E_{m,1}
\end{equation}
for all $n\, \geq\, 0$ \cite[page 259, Claim 2.7.1]{La}. Consequently,
$$
\frac{\text{degree}(E_{m,1})}{p^{m}\cdot
\text{rank}(E_{m,1})} \, \in\, {\mathbb Q}
$$
is independent of $m$ as long as $m$ is sufficiently large.
This well--defined rational number
$\frac{\text{degree}(E_{m,1})}{p^{m}\cdot
\text{rank}(E_{m,1})}$,
where $m$ is sufficiently large, is denoted by $L_{\rm max}(E)$.

A torsionfree coherent sheaf $E$ on $Y$ is called \textit{strongly
semistable} if $(F^m_Y)^*E$ is semistable for all $m\, \geq\,0$.

\begin{lemma}\label{tensor}
Let $E$ and $E'$ be vector bundles over $Y$ such that
$L_{\rm max}(E)+L_{\rm max}(E')\, <\, 0$. Then
$$
H^0(Y,\, E\otimes E') \, =\, 0\, .
$$
\end{lemma}

\begin{proof}
{}From Eqn. \eqref{isop} it follows immediately that $E_{m,1}$
is strongly semistable for all $m$
sufficiently large. Since the torsionfree part of the
tensor product of any two
strongly semistable sheaves is again strongly semistable
\cite[page 288, Theorem 3.23]{RR}, we conclude that
$$
\mu_{\text{max}}((F^m_Y)^*(E\otimes E')) \, =\,
\mu_{\text{max}}((F^m_Y)^*E) +
\mu_{\text{max}}((F^m_Y)^*E')\, =\,
p^m(L_{\rm max}(E)+L_{\rm max}(E'))\, <\, 0
$$
for all $m$ sufficiently large. Consequently,
$$
H^0(Y,\, (F^m_Y)^*(E\otimes E'))\, =\, 0
$$
for all $m$ sufficiently large. This immediately implies that
$$
H^0(Y,\, E\otimes E') \, =\, 0\, ,
$$
and the proof of the lemma is complete.
\end{proof}

The cotangent bundle of $Y$ will be denoted by $\Omega_Y$.

\begin{lemma}\label{vanish1}
Let $r$ be any integer satisfying the condition
$r\cdot{\rm degree}(L)\, >\, L_{\rm max}(\Omega_Y)$. Let
$E$ be any essentially finite vector bundle over $Y$. Then
$$
H^0(Y,\, E\otimes \Omega_Y\otimes (L^*)^{\otimes r})\,=\, 0\, ,
$$
where $L^*$ is the dual of $L$.
Moreover, if $r\cdot{\rm degree}(L)\, >\, p \cdot L_{\rm max}(\Omega_Y)$
then we have 
$$
H^0(Y,\, E\otimes F^*_Y\Omega_Y\otimes (L^*)^{\otimes r})\,=\, 0\, .
$$
\end{lemma}

\begin{proof}
Since any essentially finite vector bundle is strongly semistable
of degree zero \cite[page 37, Corollary 3.5]{No1}, 
we conclude that $L_{\rm max}(E)\, =\, 0$. 
Therefore, we have
$$
L_{\rm max}(E)+L_{\rm max}(\Omega_Y\otimes(L^*)^{\otimes r}) \,=\, 
L_{\rm max}(\Omega_Y\otimes(L^*)^{\otimes r}) \,=\,
L_{\rm max}(\Omega_Y)- r\cdot \text{degree}(L)\, <\, 0\, .
$$
Hence from Lemma \ref{tensor} it follows that
$H^0(Y,\, E\otimes \Omega_Y\otimes (L^*)^{\otimes r})\,=\,0$. 

The second part of the lemma follows similarly
using the fact that $L_{\rm max}(F^*_Y\Omega_Y)
\, =\, p\cdot L_{\rm max}(\Omega_Y)$. This completes
the proof of the lemma.
\end{proof}

\begin{lemma}\label{vanish2}
Let $D\, \subset\, Y$ be a smooth divisor in the complete
linear system $\vert L^{\otimes r}\vert\, :=\,
{\mathbb P}H^0(Y,\, L^{\otimes r})^*$. Let
${\mathcal I}\, \subset\, {\mathcal O}_Y$ be the ideal
sheaf defining $D$. Then
$$
H^0(D,\, (E\otimes({\mathcal I}/{\mathcal I}^2))\vert_D)\,=\,0
$$
and 
$$
H^0(D,\, (E\otimes F_Y^*({\mathcal I}/{\mathcal I}^2))\vert_D)\,=\,0
$$
for any essentially finite vector bundle $E$ over $Y$.
\end{lemma}

\begin{proof}
Since $E$ is essentially finite, the restriction $E\vert_D$ is an
essentially finite vector bundle over $D$.

Since $D\, \in\, \vert L^{\otimes r}\vert$, the Poincar\'e
adjunction formula says that
the line bundle over $D$ defined by ${\mathcal I}/{\mathcal I}^2$
coincides with the restriction of
the line bundle $L^{-r} \, :=\, (L^*)^{\otimes r}$ to $D$. 
As $L$ is ample, the line bundle $L^{-r}\vert_D$
over $D$ is of negative degree. Therefore, it now
follows from Lemma \ref{tensor} that
$H^0(D,\, (E\otimes({\mathcal I}/ {\mathcal I}^2))\vert_D)\,=\,0$.

Since $F^*_Y L \,=\, L^{\otimes p}$, the
restriction of the line bundle $F^*_Y L^{-r}$
to $D$ is of negative degree. Using this the second
vanishing result follows similarly. This completes
the proof of the lemma.
\end{proof}

Let $G$ and $H$ be group schemes defined over the algebraically closed
field $k$. We will denote by $G$--rep (respectively, $H$--rep) the
category of all finite dimensional left $G$--modules (respectively,
$H$--modules). Given a homomorphism of group schemes 
\begin{equation}\label{gh}
\rho_0\,:\, G \,\longrightarrow \,H
\end{equation}
we have a contravariant functor
\begin{equation}\label{ghr}
{\widetilde\rho}_0\, :\, H\mbox{-rep}\,\longrightarrow\,G\mbox{-rep}
\end{equation}
that considers a $H$--module as a $G$--module through $\rho_0$.

We now reproduce Proposition 2.21 of \cite{DM}.

\begin{proposition}[DM, page 139, Proposition 2.21]\label{inj-surj}
Let $\rho_0$ be the homomorphism in Eqn. \eqref{gh} and
${\widetilde\rho}_0$ the corresponding functor in Eqn. \eqref{ghr}.
\begin{enumerate}
\item{} The homomorphism $\rho_0$ is surjective (faithfully
flat) if and only if the following two conditions
hold: the functor ${\widetilde \rho}_0$
is fully faithful and for each exact sequence 
$$
0\, \longrightarrow\, W'
\, \longrightarrow\,{\widetilde \rho}_0(V) \, \longrightarrow\,
 W'' \, \longrightarrow\, 0
$$
of $G$--representations, where $V$ is a $H$--representation,
there is an exact sequence
$$
0\, \longrightarrow\, V'
\, \longrightarrow\,V \, \longrightarrow\,
 V'' \, \longrightarrow\, 0
$$
of $H$--representations and a commutative diagram
$$
\begin{matrix}
0 &\longrightarrow & W'
& \longrightarrow &{\widetilde \rho}_0(V) & \longrightarrow &
W''& \longrightarrow & 0\\
&&\Big\downarrow \cong && \Big\downarrow {\rm Id} &&
\Big\downarrow \cong\\
0 &\longrightarrow &{\widetilde \rho}_0(V') 
& \longrightarrow &{\widetilde \rho}_0(V) & \longrightarrow &
{\widetilde \rho}_0(V'')& \longrightarrow & 0
\end{matrix}
$$
of $G$--representations.

\item{} The homomorphism $\rho_0$ in Eqn. \eqref{gh}
is a closed immersion if and only
if for each $G$--representation $V$ there exists some
$H$--representation $W$ such that $V$ is a subquotient of
the $G$--representation ${\widetilde \rho}_0(W)$.
\end{enumerate}
\end{proposition}

Let $X$ be a smooth projective variety defined over $k$. Fix an
ample line bundle $L$ over $X$. So $L$ gives a polarization
on $X$. For any nonnegative
integer $d$, we will denote by $\vert L^{\otimes d}\vert$ the
complete linear system ${\mathbb P}H^0(X,\, L^{\otimes d})^*$.

Our aim in this section is to prove the following theorem:

\begin{theorem}\label{thm1}
Let $X$ be a smooth projective variety with $\dim X \, \geq
\, 2$ and $L$ an ample line bundle over $X$. Let $D\, \in
\, \vert L^{\otimes d}\vert$ be a smooth divisor on $X$,
where $d$ is any positive integer satisfying the inequality
$$
d\, >\, L_{\rm max}(\Omega_X)/{\rm degree}(L)\, ,
$$
and let $x_0$ be a $k$--rational point of $D$.
Then the homomorphism between fundamental group schemes
$$
\pi(D,x_0) \,\longrightarrow\, \pi(X,x_0)
$$
induced by the inclusion map $D\, \hookrightarrow\, X$
is surjective.
\end{theorem}

\begin{proof}
We will use the criterion in Proposition \ref{inj-surj}(1)
for surjectivity.

Let $E$ be an essentially finite vector bundle over $X$. From
Proposition 3.10 in \cite[page 40]{No1} we
know that $E$ is trivializable over a
principal bundle over $X$ whose structure group scheme is finite.
Taking the local group scheme corresponding to the closure of the
identity element (see the proof of Proposition 7 in
\cite[Chapter II]{No2}) we
conclude that there is a connected \'etale Galois covering
$$
f_1\, :\, X_1 \, \longrightarrow\, X
$$
such that
$f^*_1 E$ is an $F$--trivial vector bundle over $X_1$. See
Definition \ref{def1} for $F$--trivial vector bundles.

Let $V$ be another essentially finite vector bundle over $X$. Let
$$
f_2\, :\, X_2 \, \longrightarrow\, X
$$
be a connected \'etale Galois covering such that $f^*_2 V$ is
an $F$--trivial vector bundle over $X_2$. Let
\begin{equation}\label{Z}
Z\, \subset\, X_1\times_X X_2
\end{equation}
be any connected component of the fiber product. Set
\begin{equation}\label{def}
f\, :=\, f_1\times f_2 \, :\, Z \, \longrightarrow\, X\, .
\end{equation}
Therefore, both $f^*E$ and $f^* V$ are $F$--trivial vector bundles
over the connected smooth projective variety $Z$.

Let $D\, \in\,\vert L^{\otimes d}\vert$
be a smooth divisor on $X$, where
$d\, >\, L_{\rm max}(\Omega_X)/{\rm degree}(L)$.
Set $D'\, :=\, f^{-1}(D)$. Since $D\, \subset\,X$ is an ample
smooth divisor, and $\dim X\, \geq\, 2$,
the divisor $D'$ is irreducible and smooth. Thus, the restriction
\begin{equation}\label{fD}
f_D \, :=\, f\vert_{D'}\, :\, D'\, \longrightarrow\, D
\end{equation}
is also a connected \'etale Galois covering.

In our set--up,
the first condition in Proposition \ref{inj-surj}(1) says that
the restriction homomorphism
\begin{equation}\label{che}
H^0(X,\, \mathcal{H}om(E,V))\, \longrightarrow\, H^0(D,\,
\mathcal{H}om(E,V)\vert_D)
\end{equation}
is an isomorphism.

To prove this, set $E_1\, :=\, f^*E$ and $V_1\, :=\,
f^*V$, where $f$ is the morphism defined in Eqn. \eqref{def}.
Since $f$ and $f_D$ (defined Eqn. \eqref{fD})
are connected \'etale Galois coverings, we have
$$
H^0(X,\, \mathcal{H}om(E,V)) \, =\,
H^0(Z,\, \mathcal{H}om(E_1,V_1))^\Gamma
$$
and
$$
H^0(D,\, \mathcal{H}om(E\vert_D,V\vert_D)) \, =\,
H^0(D' , \,\mathcal{H}om(E_1\vert_{D'},V_1\vert_{D'}))^\Gamma\, ,
$$
where $\Gamma\, :=\, \text{Gal}(f)$ is the Galois group for
the covering $f$, which is also the Galois group for $f_D$.
Therefore, to prove that Eqn. \eqref{che} is an isomorphism
it suffices to show that the restriction homomorphism
$\mathcal{H}om(E_1,V_1)\, 
\longrightarrow\, \mathcal{H}om(E_1\vert_{D'},V_1\vert_{D'})$
gives an isomorphism
\begin{equation}\label{che2}
H^0(Z,\, \mathcal{H}om(E_1,V_1)) \, =\, H^0(D' , \,
\mathcal{H}om(E_1\vert_{D'},V_1\vert_{D'}))\, .
\end{equation}
Note that the restriction homomorphism
$$
\mathcal{H}om(E_1,V_1)\,
\longrightarrow\,\mathcal{H}om(E_1\vert_{D'},V_1\vert_{D'})
$$
is $\Gamma$--equivariant. If Eqn. \eqref{che2} holds, then taking
the $\Gamma$--invariants of both sides of Eqn. \eqref{che2}
it follows that the homomorphism in Eqn. \eqref{che} is an
isomorphism.

Let
$$
F_Z\, :\, Z\, \longrightarrow\, Z
$$
be the Frobenius morphism of the variety $Z$
in Eqn. \eqref{Z}. Since both $E_1$ and $V_1$ are $F$--trivial
vector bundles over $Z$, there is a nonnegative integer
$m$ such that both $(F^m_Z)^*E_1$ and $(F^m_Z)^*V_1$ are
trivializable vector bundles.

We proceed by induction on $m$. If $m \, =\, 0$, then
both $E_1$ and $V_1$ are trivializable vector bundles over $Z$.
In that case,
$$
H^0(Z,\, \mathcal{H}om(E_1,V_1)) \, = \,
\text{Hom}_k((E_1)_{x'_0}
\, , (V_1)_{x'_0}) \, = \, H^0(D' , \,
\mathcal{H}om(E_1\vert_{D'},V_1\vert_{D'}))\, ,
$$
where $x'_0$ is any $k$--rational point of $D'$.
Therefore, Eqn. \eqref{che2} is valid if $m\, =\, 0$.

Assume that Eqn. \eqref{che2} is valid for all $E$, $V$ and $f$
as above for which $m\, \leq\, n_0-1$.
We will show that Eqn. \eqref{che2} is valid for all $E$, $V$ and $f$
for which $m\,= n_0$. It should be emphasized that we are
not fixing $f$ in the induction hypothesis. The induction
hypothesis says that for any triple $(E\, , V\, ,f)$ as
above such that $(F^m_Z)^*E_1$ and $(F^m_Z)^*V_1$ are
trivializable vector bundles,
Eqn. \eqref{che2} is valid provided $m\, \leq\, n_0-1$.

Take any triple $E$, $V$ and $f$ as above for which $m\, =\, n_0$.
Set
$$
E' \,:=\, F^*_Z E_1
$$
and 
$$
V'\,:=\,F^*_Z V_1\, .
$$
Since $E'$ and $V'$
are pull backs under the Frobenius morphism $F_Z$,
by Cartier, there are canonical
connections $\nabla^{E'}$ and $\nabla^{V'}$ on $E'$ and $V'$
respectively, with the property that the $p$--curvature of the
connections $\nabla^{E'}$ and $\nabla^{V'}$ vanish (see
\cite[Section 5]{Ka}). Furthermore, there is a natural isomorphism 
\begin{equation}\label{e6}
\mathcal{H}om_{{\mathcal O}_Z}(E_1\, , V_1) \, =\,
\mathcal{H}om_{{\mathcal O}_Z}((E',\nabla^{E'})\, ,(V',\nabla^{V'}))
\end{equation}
between the sheaf of ${\mathcal O}_Z$--linear
homomorphisms from $E_1$ to $V_1$
and the sheaf of connection preserving ${\mathcal
O}_Z$--linear homomorphisms from
$E'$ to $V'$ \cite[page 190, Theorem 5.1.1]{Ka}.

We will use the polarization $f^*L$ on $Z$ to define
degree of coherent sheaves on $Z$. For notational
convenience, given any vector bundle
$W$ on $Z$, the vector bundle $W\otimes_{{\mathcal O}_Z}
{\mathcal O}_Z(-D')$ will be denoted by $W(-D')$.
We have
\begin{equation}\label{131}
L_{\rm max}(\mathcal{H}om(E_1,V_1(-D'))) \, =\,
L_{\rm max}(V_1) - L_{\rm max}(E_1) +
\text{degree}({\mathcal O}_Z(-D'))\, .
\end{equation}
We note that as both $E_1$ and $V_1$ are $F$--trivial,
\begin{equation}\label{132}
L_{\rm max}(V_1)\, =\, 0\, =\, L_{\rm max}(E_1)\, .
\end{equation}
Combining Eqn. \eqref{131} and Eqn. \eqref{132},
$$
L_{\rm max}(\mathcal{H}om(E_1,V_1(-D')))
\, =\, \text{degree}({\mathcal O}_Z(-D')) \, <\, 0\, ,
$$
where $D'\, =\, f^{-1}(D)$.
Hence from Lemma \ref{tensor} we conclude that
\begin{equation}\label{eq7}
H^0(Z,\, \mathcal{H}om(E_1,V_1(-D')))\,=\, 0\, .
\end{equation}

Using Eqn. \eqref{eq7} in the left exact sequence of
global sections for the exact sequence of sheaves
$$
0\,\longrightarrow\,\mathcal{H}om(E_1,V_1(-D'))
\,\longrightarrow\,\mathcal{H}om(E_1,V_1)
\,\longrightarrow\,\mathcal{H}om(E_1\vert_{D'},
V_1\vert_{D'}) \,\longrightarrow\, 0
$$
it follows that the restriction homomorphism
\begin{equation}\label{resho}
H^0(Z,\, \mathcal{H}om(E_1,V_1)) \, \longrightarrow\, H^0(D'
, \, \mathcal{H}om(E_1\vert_{D'},V_1\vert_{D'}))
\end{equation}
is injective. Hence Eqn. \eqref{che2} holds if the
homomorphism in Eqn. \eqref{resho} is surjective.

Take any homomorphism of vector bundles
$$
\phi \,\in\, H^0(D' ,\, \mathcal{H}om_{{\mathcal O}_{D'}}(
E_1\vert_{D'}, V_1\vert_{D'}))\, .
$$
over $D'$. This homomorphism $\phi$ defines a homomorphism
\begin{equation}\label{tD}
{\widetilde \phi}\, :=\,
F^*_{D'} \phi \, :\, E'\vert_{D'} \,\longrightarrow \,
V'\vert_{D'}\, ,
\end{equation}
where $F_{D'}\, :\, D'\, \longrightarrow\, D'$ is the 
Frobenius morphism of $D'$; note that $E'\vert_{D'}$
(respectively, $V'\vert_{D'}$) is identified with $F^*_{D'}
(E_1\vert_{D'})$ (respectively, $F^*_{D'} (V_1\vert_{D'})$).

Let $\nabla'$ (respectively, $\nabla''$) be the restriction
to $D'$ of the Cartier
connection $\nabla^{E'}$ (respectively, $\nabla^{V'}$)
on $E'$ (respectively, $V'$). Since Cartier
connection is compatible with pull back
operation, the connection $\nabla'$ (respectively,
$\nabla''$)
coincides with the Cartier connection on the
Frobenius pull back $F^*_{D'} (E_1\vert_{D'})\, =\, E'\vert_{D'}$
(respectively, $F^*_{D'} (V_1\vert_{D'})\, =\, V'\vert_{D'}$).
The homomorphism ${\widetilde \phi}$ in Eqn. \eqref{tD}, being
a pulled back one, intertwines the two connections $\nabla'$ and
$\nabla''$. In other words, the following identity holds:
\begin{equation}\label{e5}
({\widetilde \phi}\otimes\text{Id}_{\Omega_{D'}})\circ
\nabla'\, =\, \nabla''\circ
{\widetilde \phi}\, ,
\end{equation}
where $\Omega_{D'}$ is the cotangent bundle of $D'$.
Both sides of Eqn. \eqref{e5} are $k$--linear homomorphisms from
$E'\vert_{D'}$ to $V'\vert_{D'}\otimes \Omega_{D'}$.

Using the induction hypothesis, which says that Eqn. \eqref{che2}
is valid for all $E$, $V$ and $f$ with $m\, \leq\, n_0-1$,
the restriction homomorphism gives an isomorphism
\begin{equation}\label{eq71}
H^0(Z,\, \mathcal{H}om(E'\, ,V'))\, \stackrel{=}{\longrightarrow}\,
H^0(D',\, \mathcal{H}om(E'|_{D'},V'|_{D'}))\, .
\end{equation}
Note that as $F^*_Z E_1 \, =\, f^*F^*_X E$ and $F^*_Z V_1
\, =\, f^*F^*_X V$, and $m\,=\, n_0$ for the triple
$(E\, , V\, , f)$,
the induction hypothesis indeed
applies for the triple $(F^*_X E\, , F^*_X V\, ,f)$ giving
the isomorphism in Eqn. \eqref{eq71}.

Let
\begin{equation}\label{e7}
{\widetilde \phi}' \,:\, E' \,\longrightarrow\, V'
\end{equation}
be the homomorphism that corresponds to ${\widetilde \phi}$
(defined in Eqn. \eqref{tD}) by the isomorphism
in Eqn. \eqref{eq71}. Therefore, ${\widetilde \phi}$ is
the restriction of ${\widetilde \phi}'$ to $D'$.

To prove that Eqn. \eqref{che2} is valid for all $E$, $V$ and $f$
for which $m\, \leq\, n_0$, it suffices to show that 
the homomorphism ${\widetilde \phi}'$ in Eqn. \eqref{e7}
intertwines the 
connections $\nabla^{E'}$ and $\nabla^{V'}$. Indeed, if
${\widetilde \phi}'$ intertwines $\nabla^{E'}$ and $\nabla^{V'}$,
then using Eqn. \eqref{e6} the homomorphism ${\widetilde \phi}'$
descends to a homomorphism from $E_1$ to $V_1$. Sending any
$\phi$ to the descend of the homomorphism ${\widetilde \phi}'$
constructed from $\phi$ in Eqn. \eqref{e7}, the inverse of
the restriction homomorphism in Eqn. \eqref{resho} is obtained.

To prove that ${\widetilde \phi}'$ intertwines $\nabla^{E'}$ and 
$\nabla^{V'}$, consider
\begin{equation}\label{def.g}
\gamma\, :=\,
({\widetilde \phi}'\otimes\text{Id}_{\Omega_Z})
\nabla^{E'}-\nabla^{V'}{\widetilde
\phi}' \, \in\, H^0(Z,\, \mathcal{H}om(E', V'\otimes \Omega_Z))\, .
\end{equation}
(This homomorphism $\gamma\, :\, E' \,\longrightarrow\, V'\otimes
\Omega_Z$ is clearly ${\mathcal O}_Z$--linear.)
The homomorphism ${\widetilde \phi}'$ intertwines $\nabla^{E'}$
and $\nabla^{V'}$ if and only if $\gamma\,= \, 0$. Let
\begin{equation}\label{def.gp}
\gamma'\, \in\, H^0(D',\, (\mathcal{H}om(E',V')\otimes
\Omega_{Z})\vert_{D'}) 
\end{equation}
be the restriction to $D'$ of $\gamma$ defined in Eqn. \eqref{def.g}.

As mentioned earlier,
the degree of a coherent sheaf on $Z$ will be defined using
$f^*L$. Therefore, we have $\text{degree}(f^*W) \, =\,
\text{degree}(f)^{\dim X} \text{degree}(W)$ for any coherent
sheaf $W$ on $X$. Also,
$f^*\Omega_X \, =\, \Omega_Z$. Therefore, the given condition
in the theorem that $d\, >\, L_{\rm max}(\Omega_X)/{\rm degree}(L)$,
which is equivalent to the condition
${\rm degree}({\mathcal O}_X(D))\, >\, L_{\rm max}(\Omega_X)$,
implies that
\begin{equation}\label{eqq1}
\text{degree}({\mathcal O}_Z(D')) \, >\,
L_{\rm max}(\Omega_Z)\, .
\end{equation}
In view of Eqn. \eqref{eqq1}, from Lemma \ref{vanish1} it
follows that
$$
H^0(Z,\, \mathcal{H}om(E', V'\otimes \Omega_Z)(-D'))\, =\, 0
$$
(the vector bundle $\mathcal{H}om(E', V')$ is $F$--trivial
as both $E'$ and $V'$ are so). Therefore, considering the
left exact sequence of global sections for the exact
sequence of sheaves
$$
0\, \longrightarrow\,\mathcal{H}om(E', V'\otimes \Omega_Z)(-D')
\, \longrightarrow\,\mathcal{H}om(E', V'\otimes \Omega_Z)
\, \longrightarrow\,(\mathcal{H}om(E',V')\otimes
\Omega_{Z})\vert_{D'} \, \longrightarrow\, 0
$$
we conclude that the restriction homomorphism
\begin{equation}\label{145}
H^0(Z,\, \mathcal{H}om(E', V'\otimes \Omega_Z))\, \longrightarrow
\, H^0(D',\, (\mathcal{H}om(E',V')\otimes \Omega_{Z})\vert_{D'})
\end{equation}
is injective.

Consequently, the homomorphism $\gamma$ constructed in Eqn.
\eqref{def.g} vanishes if its restriction $\gamma'$ (defined
in Eqn. \eqref{def.gp}) vanishes.

Let ${\mathcal I}\, \subset\, {\mathcal O}_Z$ be the ideal sheaf
defining the smooth divisor $D'$.
We have an exact sequence of vector bundles over $D'$
\begin{equation}\label{eq5}
0\, \longrightarrow\,{\mathcal I}/ {\mathcal I}^2
\, \longrightarrow\, \Omega_Z\vert_{D'}\, \longrightarrow\,
\Omega_{D'} \, \longrightarrow\, 0\, .
\end{equation}
Tensoring the exact sequence in Eqn. \eqref{eq5} with
$\mathcal{H}om(E'\vert_{D'}\, ,V'\vert_{D'})$,
and then taking global sections, we obtain
the following exact sequence:
\begin{equation}\label{e4}
0\, \longrightarrow\, H^0(D',\, \mathcal{H}om(E'\vert_{D'}\,
,V'\vert_{D'})\otimes
({\mathcal I}/ {\mathcal I}^2))
\, \longrightarrow\, H^0(D',\,\mathcal{H}om(E'\vert_{D'}\,
,V'\vert_{D'})\otimes \Omega_Z\vert_{D'})
\end{equation}
$$ 
\stackrel{\beta}{\longrightarrow}\, H^0(D',\,
\mathcal{H}om(E'\vert_{D'}\, ,V'\vert_{D'})\otimes \Omega_{D'})\, .
$$
We note that the identity in
Eqn. \eqref{e5} is equivalent to the assertion
that $\beta (\gamma') \, =\, 0$, where
$\gamma'$ is the section defined in Eqn. \eqref{def.gp}
and $\beta$ is the homomorphism in Eqn. \eqref{e4}. On the
other hand, as both $E'$ and $V'$ are $F$--trivial, using the
first part of Lemma \ref{vanish2} we conclude that the term
$H^0(D',\, \mathcal{H}om(E'\vert_{D'}\, ,V'\vert_{D'})\otimes
({\mathcal I}/ {\mathcal I}^2))$ in Eqn. \eqref{e4} vanishes.
In other words, $\beta$ is injective. Since $\beta$ is injective
with $\beta (\gamma')\, =\, 0$, we conclude that $\gamma'\, =\, 0$.

We noted earlier that
$\gamma$ constructed in Eqn. \eqref{def.g} vanishes
if $\gamma'\, =\, 0$. Therefore, we have proved that
$\gamma\, =\, 0$. Since $\gamma\, =\, 0$,
the homomorphism ${\widetilde \phi}'$ in
Eqn. \eqref{e7} intertwines the
connections $\nabla^{E'}$ and $\nabla^{V'}$.
We noted earlier that this implies that 
Eqn. \eqref{che2} is valid for all $E$, $V$ and $f$
for which $m\, \leq\, n_0$.

Using induction on $m$ it now follows
that Eqn. \eqref{che2} is valid.
It was also shown earlier that 
Eqn. \eqref{che2} implies that the homomorphism in
Eqn. \eqref{che} is an isomorphism. This completes
the proof of the assertion that the homomorphism
\begin{equation}\label{133}
\pi(D,x_0) \,\longrightarrow\, \pi(X,x_0)
\end{equation}
in the statement of the theorem satisfies the first condition
in Proposition \ref{inj-surj}(1).

We will now prove that the homomorphism
in Eqn. \eqref{133} also satisfies
the second condition in Proposition \ref{inj-surj}(1).
For that we will use the above set--up, and we will again
employ induction on $m$.

Let $V$ be an essentially finite vector bundle over $X$
such that $(F^m_Z)^*f^*V$ is a trivial vector bundle over $Z$,
where $f$ and $Z$ are as above. Let
\begin{equation}\label{e0}
E_0\, \subset\, V\vert_D
\end{equation}
be an essentially finite subbundle over $D$. Therefore,
the vector bundle $(F^m_{D'})^*f^*_D E_0$ over $D'\, :=\,
f^{-1}(D)$, where $f_D$ is the restriction
morphism in Eqn. \eqref{fD} and $F_{D'}$ is the
Frobenius morphism of $D'$ (as in Eqn. \eqref{tD}),
is an essentially finite subbundle of the trivializable
vector bundle $((F^m_Z)^*f^*V)\vert_{D'}$. Since any
subbundle of degree zero of a trivializable vector bundle
is also trivializable, and any essentially finite vector
bundle is of degree zero, we conclude that $(F^m_{D'})^*f^*_D E_0$
is a trivializable vector bundle over $D'$.

If $m\, =\, 0$, then it is easy to see that $E_0$ is the restriction
to $D$ of an essentially finite subbundle of $V$ over $X$. Indeed,
this follows immediately from the fact that the natural
homomorphism of \'etale fundamental groups
$$
\pi_1(D,x_0) \,\longrightarrow\, \pi_1(X,x_0)
$$
is surjective (a proof of this is given in Section \ref{sec.ex}).

Assume that for all triples of the form $(V\, , E_0\, , f)$ as
above with
$m\, < \, n_0$, the vector bundle $E_0$ is the restriction
to $D$ of an essentially finite subbundle of $V$. We
will show that for all triples $(V\, , E_0\, , f)$ as above
with $m\, = \, n_0$, the vector bundle $E_0$ is the restriction
to $D$ of an essentially finite subbundle of $V$.

Take any triple $(V\, , E_0\, , f)$ as above for which
$m\, = \, n_0$. Let
\begin{equation}\label{134}
V'\, :=\, F^*_{Z}f^* V
\end{equation}
be the vector bundle over $Z$ and
\begin{equation}\label{104}
E'_0\, :=\, F^*_{D'}f^*_D E_0
\end{equation}
the vector bundle over $D'$. Since $f\circ F_Z \, =\,
F_X\circ f$, we have
$$
(F^m_Z)^*f^*V\, =\, (F^{m-1}_Z)^*f^*F^*_X V\, .
$$
Also, since $F_X\circ \iota_D\, =\, \iota_D\circ F_D$,
where
$$
\iota_D\, :\, D\, \hookrightarrow\, X
$$
is the inclusion morphism and $F_D\, :\, D\, \longrightarrow\, D$
is the Frobenius morphism, the inclusion
$$
E_0\, \subset\, V\vert_D
$$
in Eqn. \eqref{e0} induces an injective homomorphism of vector
bundles
$$
F^*_D E_0\, \subset\, (F^*_XV)\vert_D\, .
$$
Consequently, the induction hypothesis applies for the
triple $(F^*_X V\, , F^*_{D} E_0\, , f)$. (Note that
$F^*_{D} E_0$ is essentially finite as $E_0$ is so.)
Therefore, there is an essentially finite subbundle
$$
E_1\, \subset\, F^*_X V
$$
such that $E_1\vert_D \, =\, F^*_{D}E_0\, \subset\,
(F^*_X V)\vert_D$.

Set
\begin{equation}\label{ep}
E'\, :=\, f^* E_1\, \subset\, V'\, .
\end{equation}
Therefore, the two subbundles of $V'\vert_{D'}$, namely
$E'\vert_{D'}$ and $E'_0$, coincide, where $E'_0$ is
defined in Eqn. \eqref{104}. Note that
\begin{equation}\label{105}
E'_0\, =\, F^*_{D'}f^*_D E_0
\, \subset\, F^*_{D'}f^*_D(V\vert_D)\, =\,
F^*_{D'}((f^*V)\vert_{D'})\, =\,
(F^*_Zf^*V)\vert_{D'}\, =\, V'\vert_{D'}\, .
\end{equation}

As before, let $\nabla^{V'}$ denote the Cartier connection on
the vector bundle $V'$ defined in Eqn. \eqref{134}. Let
\begin{equation}\label{ded}
\delta\, :\, E'\, \hookrightarrow\, V'\,
\stackrel{\nabla^{V'}}{\longrightarrow}\, V'\otimes \Omega_Z
\, \longrightarrow\, (V'/E')\otimes \Omega_Z
\end{equation}
be the composition homomorphism giving the
second fundamental form of the subbundle $E'$ (defined in
Eqn. \eqref{ep}) for the Cartier connection
$\nabla^{V'}$. From Leibniz identity, this composition
homomorphism is ${\mathcal O}_Z$--linear. In other words,
\begin{equation}\label{ded2}
\delta\, \in\, H^0(Z,\, \mathcal{H}om(E'\, ,(V'/E')\otimes
\Omega_Z))\, .
\end{equation}
This homomorphism $\delta$ vanishes if and only if the
subbundle $E'$ is the pull back, by $F_Z$, of a subbundle
of $f^* V$ \cite[page 190, Theorem 5.1]{Ka}.

Since the vector bundle $V'$ is $F$--trivial and
$E'$ is an essentially finite subbundle of $V'$, it follows
immediately that both $E'$ and $V'/E'$ are $F$--trivial
vector bundles (a subbundle of degree zero of a
trivial vector bundle is trivializable and the quotient is also
trivializable). As we noted earlier, the inequality in Eqn.
\eqref{eqq1} follows from the given condition on $d$.
Since $E'$ and $V'/E'$ are $F$--trivial,
using Lemma \ref{vanish1} it follows that
$$
H^0(Z,\, \mathcal{H}om(E'\, ,(V'/E')\otimes
\Omega_Z)(-D'))\, =\, 0\, .
$$
In view of this, considering the left exact sequence
of global sections for the exact sequence
$$
0\,\longrightarrow\,\mathcal{H}om(E'\, ,(V'/E')\otimes
\Omega_Z)(-D')\,\longrightarrow\, \mathcal{H}om(E'\, ,(V'/E')
\otimes\Omega_Z)
$$
$$
\,\longrightarrow\, \mathcal{H}om(E'\, ,
(V'/E')\otimes \Omega_Z)\vert_{D'} \,\longrightarrow\,0
$$
we conclude that the restriction homomorphism
\begin{equation}\label{135}
\psi\, :\, H^0(Z,\, \mathcal{H}om(E'\, ,(V'/E')\otimes
\Omega_Z))\, \longrightarrow\, H^0(D', \, \mathcal{H}om(E'\, ,
(V'/E')\otimes \Omega_Z)\vert_{D'})
\end{equation}
is injective.

Set
\begin{equation}\label{dep}
\delta'\, :=\, \psi(\delta)\, \in\, H^0(D', \,
\mathcal{H}om(E'\, ,(V'/E')\otimes\Omega_Z)\vert_{D'})\, ,
\end{equation}
where $\psi$ is the homomorphism in Eqn. \eqref{135}
and $\delta$ is the section in Eqn. \eqref{ded2}.
In other words, $\delta'$ is the restriction of $\delta$
to $D'$.

Since the homomorphism $\psi$ in Eqn. \eqref{135} is injective, the
section $\delta$ vanishes if
$\delta'$, defined in Eqn. \eqref{dep}, vanishes.

Tensoring the short exact sequence in Eqn. \eqref{eq5} with
$\mathcal{H}om(E'\vert_{D'}\, ,(V'/E')\vert_{D'})$, and then taking
global sections, we get a left exact sequence
\begin{equation}\label{e8}
H^0(D',\, \mathcal{H}om(E'\vert_{D'}\,
,(V'/E')\vert_{D'})\otimes ({\mathcal I}/ {\mathcal I}^2))
\, \longrightarrow\, H^0(D',\, \mathcal{H}om(E'\vert_{D'}\,
,(V'/E')\vert_{D'})\otimes \Omega_Z\vert_{D'})
\end{equation}
$$
\stackrel{\alpha}{\longrightarrow}\, H^0(D',\,
\mathcal{H}om(E'\vert_{D'}\, ,(V'/E')\vert_{D'})\otimes \Omega_{D'})\, .
$$
As $E'$ and $V'/E'$ are $F$--trivial vector bundles, the vector
bundles $E'\vert_{D'}$ and $(V'/E')\vert_{D'}$ are also $F$--trivial.
Hence the first part of Lemma \ref{vanish2} gives that
$$
H^0(D',\, \mathcal{H}om(E'\vert_{D'}\,
,(V'/E')\vert_{D'})\otimes ({\mathcal I}/ {\mathcal I}^2))
\, =\, 0\, .
$$
Therefore, the homomorphism $\alpha$ in Eqn. \eqref{e8}
is injective.

Recall that $E'\vert_{D'}\, =\, E'_0$, or in other words,
the restriction $E'\vert_{D'}$
is the pull back, by the Frobenius morphism $F_{D'}$,
of the subbundle $f^*_D E_0 \, \subset\, (f^*V)\vert_{D'}$.
Therefore, the second fundamental form of the subbundle
$$
E'\vert_{D'}\, \subset\, V'\vert_{D'}\,=\, F^*_{D'}((f^*V)\vert_{D'})
$$
(see Eqn. \eqref{105})
for the Cartier connection on $F^*_{D'}f^*_D (V\vert_D)$
vanishes identically. This immediately implies
that $\alpha(\delta')\, =\, 0$, where $\delta'$ is defined
in Eqn. \eqref{dep} and $\alpha$ is the homomorphism in
Eqn. \eqref{e8}.

Since $\alpha$ is injective (this was proved above), and
$\alpha(\delta')\, =\, 0$, we conclude that $\delta'\, =\, 0$.

We saw earlier that if $\delta'\, =\, 0$, then
the homomorphism $\delta$ constructed in Eqn. \eqref{ded} vanishes.
Therefore, we have proved that $\delta\, =\, 0$. Since
$\delta\, =\, 0$,
the subbundle $E'\, \subset\, V'$ is the pull back,
by $F_Z$, of a subbundle of $f^* V$. Let
\begin{equation}\label{epp}
E''\, \subset\, f^* V
\end{equation}
be the subbundle such that $F^*_Z E'' \, =\, E'\, \subset\, V'$.
Since $E'$ is essentially finite, it follows that
the vector bundle $E''$ is also essentially finite. Indeed. as
$E'$ is essentially finite, there is an \'etale Galois cover
$Z'$ of $Z$ such that the pull back of $E'$ to $Z'$
is $F$--trivial, hence the pull back $E''$ to
$Z'$ is $F$--trivial, which, using Proposition
\ref{pr.e.f.}), implies that $E''$ is essentially finite.

We will show that the natural
action of the Galois group $\Gamma\, :=\,
\text{Gal}(f)$ on $f^* V$ leaves the subbundle $E''$ invariant. For
this, take any element $g\,\in\, \Gamma$ of the Galois group. Let
\begin{equation}\label{Hg}
H_g\, :\, E''\oplus (g^{-1})^*E'' \, \longrightarrow\, f^*V
\end{equation}
be the homomorphism defined by $H_g(v,w) \, =\, \iota(v) +
g(((g^{-1})^*\iota)(w))$, where $\iota\, :\, E''\, \hookrightarrow\,
f^*V$ is the inclusion map in Eqn. \eqref{epp},
$$
(g^{-1})^*\iota\, :\,
(g^{-1})^*E''\, \hookrightarrow\, (g^{-1})^*f^*V \,=\, f^*V
$$
is the pull back of $\iota$ and $g(((g^{-1})^*\iota)(w))$ is the
image of $((g^{-1})^*\iota)(w)\,\in\, f^*V$ for the action of $g$
on $f^*V$. The vector bundle $E''\oplus (g^{-1})^*E''$ is
essentially finite as $E''$ is so. Therefore,
the image $H_g(E''\oplus (g^{-1})^*E'')$ is a subbundle of the
essentially finite vector bundle
$f^*V$ \cite[page 38, Proposition 3.7(b)]{No1}.

Since
$$
F^*_{D'}(E''\vert_{D'})\, =\,
(F^*_{Z}E'')\vert_{D'}\,=\, E'\vert_{D'}\,=\, E'_0\,=\,
F^*_{D'}f^*_D E_0\, ,
$$
we know that $f^*_D E_0 \, =\, E''\vert_{D'}$.
Therefore, the
action of $\Gamma$ on $(f^*V)\vert_{D'}$ leaves the subbundle
$$
f^*_D E_0 \, =\, E''\vert_{D'} \, \subset\, (f^*V)\vert_{D'}
$$
invariant. Hence we have
\begin{equation}\label{rank}
\text{rank}(H_g(E''\oplus E'')) \, =\, \text{rank}(E'')
\end{equation}
for the homomorphism $H_g$ in Eqn. \eqref{Hg}. From
Eqn. \eqref{rank} it follows immediately that the action of
$g\,\in\, \Gamma$ on $f^*V$ leaves the subbundle $E''$ invariant.

Since the subbundle $E''\, \subset\, f^*V$ is
left invariant by $\Gamma$, it 
descends as a subbundle of $V$, and furthermore, the descended
subbundle extends to $X$
the subbundle $E_0\, \subset\, V\vert_D$ over $D$. Since $E''$ is
essentially finite, it follows that the descend of $E''$ is an
essentially finite vector bundle over $X$. Indeed, as $f$ is \'etale
Galois, and $E''$ is essentially finite, there is an
\'etale Galois cover $X'$ of $X$ such that the pull back
to $X'$ of the descend of $E''$ to $X$ is $F$--trivial,
which, using Proposition \ref{pr.e.f.}, implies that the descend
of $E''$ to $X$ is essentially finite.

Thus, if $V$ is an essentially finite vector bundle over $X$ and
$E_0\,\subset\, V\vert_D$ an essentially finite
subbundle over $D$, where $D\, \in\,\vert L^{\otimes d}\vert$
is a smooth divisor on $X$ with
$d\, >\, L_{\rm max}(\Omega_X)/{\rm degree}(L)$,
then $E_0$ extends to $X$ as an
essentially finite subbundle of $V$. Therefore, the homomorphism
in Eqn. \eqref{133} satisfies
the second condition in Proposition \ref{inj-surj}(1).
Hence this homomorphism is surjective.
This completes the proof of the theorem.
\end{proof}

Combining Remark \ref{rem1} and Theorem \ref{thm1} we have the
following corollary:

\begin{corollary}\label{cor.va}
Let $X$ be a smooth projective variety of dimension at least
two and $L$ an ample line bundle over $X$.
Let $D\, \in\,
\vert L^{\otimes d}\vert$ be a smooth divisor on $X$,
where $d$ is any positive integer satisfying the inequality
$$
d\, >\, L_{\rm max}(\Omega_X)/{\rm degree}(L)\, .
$$
Then $H^1(X,\, {\mathcal O}_X(-D))\, =\, 0$.
\end{corollary}

\section{Some vanishing results}\label{sec.s.v.r.}

In this section we will prove some vanishing results which
will be used in Section \ref{sec.-inj.}.

Let $k$ be an algebraically closed field of
characteristic $p$, with $p\, >\, 0$.
Let $X$ be a smooth projective variety
defined over $k$ and $F_X$ the Frobenius morphism
of $X$. Given any vector bundle $E$ over $X$, we will construct
a natural filtration of coherent subsheaves of $F^*_XF_{X*}E$.

Let $W$ be a vector bundle over $X$ equipped with a connection
$\nabla$. Let $W_1\, \subset\, W$ be a coherent subsheaf. Set
$W_2$ to be the kernel of the second fundamental form
$$
S(W_1) \, :\, W_1\, \longrightarrow\, (W/W_1)\otimes \Omega_X
$$
of $W_1$ for the connection $\nabla$, where $\Omega_X$ is the
cotangent bundle of $X$. We recall that
$$
S(W_1)\,=\, (q_{W_1}\otimes\text{Id}_{\Omega_X})\circ\nabla\circ
\iota_{W_1}\, ,
$$
where $\iota_{W_1}\, :\, W_1 \,\hookrightarrow\, W$ is the inclusion
map and $q_{W_1}\, :\, W \, \longrightarrow\, W/W_1$ is the
quotient map. Next set $W_3$ to be the kernel of the
second fundamental form
$$
S(W_2) \, :\, W_2\, \longrightarrow\, (W/W_2)\otimes \Omega_X
$$
of $W_2$ for the connection $\nabla$. This way, for
$i\, \geq\, 1$, define
$W_{i+1}$ inductively to be the kernel of the second fundamental
form $S(W_i)$ of $W_i$ for the connection $\nabla$.
For notational convenience, set $W_0\, :=\, W$. From the construction
of this filtration
\begin{equation}\label{41}
W\, =\, W_0 \, \supset\, W_1 \, \supset\, W_2 \, \supset\, \cdots
\end{equation}
it follows immediately that
$$
\text{Image}(S(W_i)) \, \subset\, (W_{i-1}/W_i)\otimes\Omega_X
\, \subset\, (W/W_i)\otimes\Omega_X
$$
for each $i\, \geq\, 1$. Consequently, the second
fundamental form $S(W_i)$ induces a homomorphism
\begin{equation}\label{va.e1}
{\widetilde S}(W_i)\, :\, W_i/W_{i+1}\, \longrightarrow\,
(W_{i-1}/W_{i})\otimes\Omega_X
\end{equation}
for each $i\, \geq\, 1$.

Iterating the homomorphisms in Eqn. \eqref{va.e1} we have
a homomorphism
$$
{\widetilde S}(W_1)\circ\cdots\circ
{\widetilde S}(W_{i-1})\circ{\widetilde S}(W_i)\, :\,
W_i/W_{i+1}\, \longrightarrow\, (W/W_1)\otimes
\Omega^{\otimes i}_X
$$
for all $i\, \geq\, 1$. Let
\begin{equation}\label{va.e2}
S_i \, :\, W_i/W_{i+1}\, \longrightarrow\, (W/W_1)\otimes
\Omega^{\otimes i}_X
\end{equation}
be the composition homomorphism.

Take any vector bundle $E$ over $X$.
The vector bundle $F^*_XF_{X*}E$ is equipped with the Cartier
connection. First set $E_0\, :=\, F^*_XF_{X*}E$. There is a natural
surjective homomorphism
\begin{equation}\label{va.m1}
E_0\, :=\, F^*_XF_{X*}E\, \longrightarrow\, E\, .
\end{equation}
Let $E_1\, \subset\, E_0$ denote the kernel of the
homomorphism in Eqn. \eqref{va.m1}. Using the
construction of the filtration in Eqn. \eqref{41},
a filtration of $F^*_XF_{X*}E$ is obtained from the
subbundle $E_1$
and the Cartier connection on $F^*_XF_{X*}E$. Thus,
we get a filtration of coherent subsheaves
\begin{equation}\label{va.e4}
F^*_XF_{X*}E\, =:\, E_0\, \supset\, E_1 \, \supset\,
E_2 \, \supset\, \cdots\, .
\end{equation}
Since the homomorphism in Eqn. \eqref{va.m1}
is surjective, we have $E_0/E_1\, =\, E$.

\begin{proposition}\label{canfil}
Let $p$ denote the characteristic of the field $k$, and
$d\, =\, \dim X$.
With the above notation we have the following:
\begin{enumerate}
\item{} Each $E_i$ in Eqn. \eqref{va.e4} is a vector bundle over
$X$, and $E_i \, =\, 0$ if $i\, >\, (p-1)d$.

\item{} For each $i\,\geq\, 1$, the homomorphism
\begin{equation}\label{42}
S_i \, :\, E_i/E_{i+1}\,\longrightarrow\, E\otimes\Omega^{\otimes i}_X
\end{equation}
constructed as in Eqn. \eqref{va.e2} is an injective homomorphism
of vector bundles.

\item{} Set $E\, =\, {\mathcal O}_X$. Then for each $i\, \geq\, 1$,
consider the subbundle
$$
{\mathcal G}_i \subset \Omega^{\otimes i}_X
$$
given by the image of the homomorphism $S_i$ in
Eqn. \eqref{42} (for
$E\, =\, {\mathcal O}_X$).
Let $E$ be an arbitrary vector bundle over $X$.
Then the image $S_i(E_i/E_{i+1})$ coincides with the
subbundle
$$
E \otimes {\mathcal G}_i\, \subset\, E\otimes\Omega^{\otimes i}_X\, ,
$$
where ${\mathcal G}_i$ is the above subbundle of
$\Omega^{\otimes i}_X$.
\end{enumerate}
\end{proposition}

\begin{proof}
If $f\, :\, X_1\, \longrightarrow\, X_2$ is a morphism, then
$f\circ F_{X_1}\, =\, F_{X_2}\circ f$, where $F_{X_j}$, $j\,=\,1,2$,
is the Frobenius morphism of $X_j$. Furthermore, if
$V$ is a vector bundle over $X_2$, then the isomorphism
$F^*_{X_1}f^*V\, =\, f^*F^*_{X_2}V$ takes the Cartier connection
on $F^*_{X_1}f^*V$ to the pullback, by $f$, of the
Cartier connection on $F^*_{X_2}V$. Using these and the fact that
the proposition is local in nature we conclude that
it is enough to prove the proposition assuming that
$X\, =\, \text{Spec}(A)$, where $A\,:=\,k[[x_1,\ldots, x_d]]$.
We will assume that $X\, =\, \text{Spec}(A)$ unless specified
otherwise. Therefore, we
may also assume that the vector bundle $E$ over $X$ is
trivial. We will identify $E$ with $A^{\oplus r}$ equipped with
the standard basis $\{e_1\, , \ldots\, , e_r\}$. 

Therefore, the cotangent bundle $\Omega_X$ is a free $A$--module with
a basis $\{{\rm d}x_1, \cdots , {\rm d}x_d\}$.

Set $B\, :=\, k[[x_1^p,\ldots ,x_d^p]]\, \subset\, A$, which
is isomorphic to
$A$. The Frobenius morphism $F_X$ is seen as the inclusion map
$$
B \,\hookrightarrow\, A\, .
$$
Using this inclusion we see that $A$ is a free $B$--module with a 
basis given by
monomials of the form $x_1^{m_1}\cdots x_d^{m_d}$, where 
$0\, \leq\, m_j \,\leq\, p-1$.
The total degree of such a monomial $x_1^{m_1}\cdots x_d^{m_d}$
is defined to be $\sum_{j=1}^d m_j$.

We will identify $F^*_XF_{X*}E$ with $\bigoplus_{i=1}^r
(A\otimes_B A)e_i$
using the trivialization of $E$ and the identification of $F_X$
as the inclusion of $B$ in $A$. The $B$--module structure of
$$
F_{X*}E\, =\, \bigoplus_{i=1}^r Ae_i
$$
is defined by the following rule:
$$
x^p_j(\sum_{i=1}^r a_ie_i)\,=\, \sum_{i=1}^r x_j^p a_ie_i\, ,
$$
where $j\, \in\, [1\, ,d]$ and $a_i\,\in\, A$.
We note that $F_{X*} E$
is a free $B$--module with a basis
given by $\{x_1^{m_1}x_2^{m_2}\cdots x_d^{m_d}e_i\}$,
where $1\,\leq\, i\,\leq\, r$
and $0\,\leq\, m_j\, <\, p$ for each $j\, \in\, [1\, ,d]$.

The natural projection $F^*_XF_{X*} E\,\longrightarrow\, E$
in Eqn. \eqref{va.m1} is now identified with a homomorphism
$$
\bigoplus_{i=1}^r (A\otimes_B A)e_i \,\longrightarrow\,
\bigoplus_{i=1}^r Ae_i
$$
which is defined as follows:
\begin{equation}\label{va.s0}
c\otimes x^{m_1}_1\cdots x^{m_d}_d e_i\, \longmapsto\,
cx^{m_1}_1\cdots x^{m_d}_d \otimes{1} e_i\, ,
\end{equation}
where $C^i_{m_1,\ldots, m_d}\, \in\, A$.

For $i\, \in\, [1\, ,d]$, let $y_i$ be the variable defined
by $y_i\, :=\, x_i\otimes 1 -1\otimes x_i\, \in\,
A\bigotimes_B A$.

The homomorphism defined in Eqn. \eqref{va.s0} is surjective and
its kernel is generated
by linear combinations of all $e_i$ with coefficients
that are homogeneous polynomials in the
variables $y_i$ of total degree at least one.
This already shows that the subsheaf $E_1$
in Eqn. \eqref{va.e4} is a subbundle of $E_0$.

The Cartier connection operator on $F^*_XF_{X*} E$
$$
F^*_XF_{X*} E \, \longrightarrow\, (F^*_XF_{X*} E)
\otimes \Omega_X
$$
coincides with the map
$$
\bigoplus_{i=1}^r (A\otimes_B A)e_i\, \longrightarrow\,
\bigoplus_{i=1}^r (A\otimes_B A)\otimes_A \Omega_X e_i
$$
defined by
$$
\sum_{i=1}^r C^i(y_1,\ldots,y_d)e_i\,\longmapsto\,
\sum_{i=1}^r {\rm d}(C^i(y_1,\ldots,y_d))\otimes e_i\, ,
$$
where 
$$
{\rm d}(C^i(y_1,\ldots,y_d))\,:=\, \sum_{j=1}^d
\frac{\partial C^i(y_1,\ldots,y_d)}
{\partial x_j}\otimes {\rm d}x_j
$$
with $C^i(y_1,\ldots,y_d)\, \in\, A\otimes_B A$,
and
$$
\frac{\partial}{\partial x_j}\, :\, A\otimes_B A
\,\longrightarrow\, A\otimes_B A
$$
being the pullback of the $B$--derivation
$$
\frac{\partial}{\partial x_j}\, :\, A
\,\longrightarrow\, A\, .
$$

The second fundamental form of $E_1$ for the Cartier
connection
$$
S(E_1) \, :\, E_1\, \longrightarrow\, E\otimes\Omega_X
$$
is identified with the $A$--linear homomorphism
defined as follows:
$$
\sum_{i=1}^r\sum_{\{m_j\}_{j=1}^d \in [0,p-1]^d\setminus 0^d}
C_{m_1,\cdots, m_d}\otimes
x_1^{m_1}\cdots x_d^{m_d}e_i\,\longmapsto\,
\sum_{i=1}^r \sum_{j=1}^d C_{\underline{j}}e_i\otimes{\rm d}x_j\, ,
$$
where $\underline{j}\, =\, (0, \cdots , 1, \cdots ,0)
\, \in\, {\mathbb Z}^d$, i.e., $0$ everywhere except
at the $j$-th position where it is $1$. From
this expression
of the second fundamental form $S(E_1)$ it follows immediately
that the kernel of $S(E_1)$ is a free $A$--module.

For the general case, using induction we see that
for any $\ell\, \geq\, 1$, the subsheaf $E_\ell
\, \subset\, F^*_X F_{X*} E$ is
a free $A$--module generated by monomials of the form 
$y_1^{m_1}\cdots y_d^{m_d}e_i$, where $1 \le i\le r$,
$\sum_{j=1}^d m_j \,\geq\, \ell$ and $0 \leq m_j <p$ for each
$j\, \in\, [1\, , d]$. The second fundamental form
$$
S(E_\ell) \, :\, E_\ell\, \longrightarrow\, (E_{\ell -1}/E_\ell)
\otimes\Omega_X
$$
coincides with the homomorphism that sends any element
$$
\sum_{i=1}^r C^i(y_1,\ldots,y_d)e_i\, \in\, E_\ell
$$
to the image of $\sum_{i=1}^r{\rm d}(C^i(y_1,\ldots,y_d))
e_i$ in $(E_{\ell -1}/E_\ell)\otimes\Omega_X$.
Note that the
kernel of this homomorphism is a free $A$--module generated by the
monomials of the form $y_1^{m_1}\cdots y_d^{m_d}e_i$, where
$1 \le i\le r$, $\sum_{j =1}^d m_j \,\geq\, \ell+1$ and
$0\leq m_j < p$ for all $j\, \in\, [1\, ,d]$.

{}From this description
it follows that each $E_\ell$ is a vector bundle. 
It also follows that $E_\ell\, =\, 0$ if $\ell\, >\, (p-1)d$.
Therefore, the proofs of statement (1) and statement (2)
in the proposition are complete.

To complete the proof of the proposition we need to show that
the two
subbundles of $E\otimes \Omega^{\otimes\ell}_X$, namely
$\text{image}(S_\ell)$ and $E\otimes {\mathcal G}_\ell$,
coincide.

To prove this we may assume that
$X\, =\, \text{Spec}(A)$, where $A\,:=\,k[[x_1,\ldots, x_d]]$,
and we may also assume that $E\,=\, {\mathcal O}^{\oplus r}_X
\, =\, A^{\oplus r}$.
If $E\,=\, {\mathcal O}_X\, =\, A$, then 
$$
\text{image}(S_\ell) \, =\, {\mathcal G}_\ell
\, =\, E\otimes {\mathcal G}_\ell
\, \subset\, E\otimes \Omega^{\otimes\ell}_X
\, =\, \Omega^{\otimes\ell}_X
$$
for all $\ell\, \geq\, p$. Therefore,
$$
\text{image}(S_\ell)\, =\, E\otimes {\mathcal G}_\ell
\, \subset\, E\otimes \Omega^{\otimes\ell}_X
$$
if $E\,=\, {\mathcal O}^{\oplus r}_X$ and $\ell\, \geq\, p$.
This completes the proof of the proposition.
\end{proof}

The following proposition is proved using
Proposition \ref{canfil}(1).

\begin{proposition}\label{va.}
Let $X$ be a smooth projective variety of dimension $d$
defined over $k$.
Let $M$ be a nonnegative rational number such that
$$
M\, \geq\, L_{\rm max}(\Omega^{\otimes i}_X)
$$
for all $i\, \in\, [1\, , (p-1)d]$,
where $L_{\rm max}$ is defined in Section \ref{se2} and
$p$ is the characteristic of $k$. Then
for each positive integer $n$ the following inequality holds:
$$
L_{\rm max}((F^n_{X})_* {\mathcal O}_X)\, \leq \, M\, .
$$
\end{proposition}

\begin{proof}
For notational convenience, the vector bundle
$(F^{n-1}_{X})_* {\mathcal O}_X$ over $X$ will be denoted by
$W$. Consider the filtration of quotients of $F^*_XF_{X*}W$
$$
F^*_XF_{X*}W\, =:\, W_0\, \supset\, W_1 \, \supset\, W_2
\, \supset\, \cdots
$$
constructed as in Eqn. \eqref{va.e4}. Pulling it back by the
morphism $(F^{n-1}_X)^*$ we obtain a filtration of subbundles
\begin{equation}\label{va.fip}
(F^{n-1}_X)^* F^*_XF_{X*}W\, =:\, W'_0\, \supset\, W'_1 \,
\supset\, W'_2 \, \supset\, \cdots\, ,
\end{equation}
where $W'_i \, :=\, (F^{n-1}_X)^*W_i$ for all $i\, \geq\, 0$.
Since
$$
(F^{n-1}_X)^* F^*_XF_{X*}W\, =\, (F^{n}_X)^*(F^n_{X})_*
{\mathcal O}_X\, ,
$$
the filtration in Eqn. \eqref{va.fip} gives a filtration
of subbundles of $(F^n_X)^*(F^n_{X})_* {\mathcal O}_X$.

Using Proposition \ref{canfil}(1) we know that
$$
W'_i/W'_{i+1} \, =\, (F^{n-1}_X)^* (W_i/W_{i+1})\, \subset\,
(F^{n-1}_X)^*(W\otimes \Omega^{\otimes i}_X)
$$
and $W'_i/W'_{i+1}\, =\, 0$ if $i\, >\, (p-1)d$, where
$d\, =\, \dim X$. On the other hand, given any two vector
bundles $V_1$ and $V_2$ over $X$, we have
$$
L_{\text{max}}(V_1\otimes V_2) \, =\, L_{\text{max}}(V_1)+
L_{\text{max}}(V_2)
$$
(this follows from \cite[page 288, Theorem 3.23]{RR}).
Combining these we conclude that there exists a nonnegative integer
$i\, \leq\, (p-1)d$ such that
\begin{equation}\label{va.ind.}
L_{\text{max}}((F^{n}_X)^*(F^n_{X})_* {\mathcal O}_X)\, \leq\,
L_{\text{max}}((F^{n-1}_X)^*W) + L_{\text{max}}((F^{n-1}_X)^*
\Omega^{\otimes i}_X\, .
\end{equation}
(If $0\, =\, V_0\, \subset\, V_1\, \subset\,\cdots \,\subset
\, V_{m-1} \, \subset\, V_m$ is a filtration of subbundles
and $V'\, \subset\, V_m$ is a coherent subsheaf, then there is a
nonzero homomorphism of $V'$ to some quotient $V_i/V_{i-1}$.)

Let $M$ be any nonnegative rational number such that
$$
M\, \geq\, L_{\text{max}}(\Omega^{\otimes j}_X)
$$
for all $j\, \in\, [1\, , (p-1)d]$. From the definition of
$L_{\rm max}$ we have
\begin{equation}\label{155}
L_{\text{max}}((F^{n-1}_X)^*
\Omega^{\otimes j}_X)\, =\, p^{n-1} L_{\text{max}}(
\Omega^{\otimes j}_X)\,\leq\, p^{n-1}M
\end{equation}
for all $j\, \in\, [1\, , (p-1)d]$.

As $W\, =\, (F^{n-1}_{X})_* {\mathcal O}_X$, from
the inequalities in Eqn. \eqref{va.ind.} and Eqn.
\eqref{155} we have
$$
L_{\text{max}}((F^{n}_X)^*(F^n_{X})_* {\mathcal O}_X)
\, \leq\,\sum_{j=0}^{n-1} p^j M \, =\, M(p^n-1)/(p-1)\, .
$$
Since $L_{\text{max}}((F^{n}_X)^*(F^n_{X})_* {\mathcal O}_X)
\,=\, p^n L_{\text{max}}((F^n_{X})_* {\mathcal O}_X)$,
this inequality implies that
$$
L_{\text{max}}((F^n_{X})_* {\mathcal O}_X)\, \leq\,
M(p^n-1)/p^n(p-1) \, \leq\, M\, .
$$
This completes the proof of the proposition.
\end{proof}

Let
\begin{equation}\label{va.f}
f\, : Z\, \longrightarrow\, X
\end{equation}
be a connected Galois \'etale cover of the smooth
projective variety $X$. As before, the
Frobenius morphism of $Z$ will be denoted by $F_Z$.

\begin{lemma}\label{basechange}
The following diagram is Cartesian
$$
\begin{array}{rrl}
Z & \stackrel{F_Z}{\longrightarrow} & Z \\
f\Big\downarrow& &\Big\downarrow f \\
X & \stackrel{F_X}{\longrightarrow} &X
\end{array}
$$
\end{lemma}

\begin{proof}
As the above diagram is commutative, we obtain a morphism
$$
g\, :\, Z\,\longrightarrow\, X\times_X Z
$$
to the fiber product. The composition
$$
Z\, \stackrel{g}{\longrightarrow}\, X\times_X Z
\, \stackrel{\text{pr}_X}{\longrightarrow}\, X
$$
coincides with $f$, and hence it is
an \'etale morphism. On the other hand, the projection
$\text{pr}_X$ is \'etale as $f$ is so. Therefore,
the morphism $g$ is \'etale \cite[page 24, Corollary 3.6]{Mi}.
Since the \'etale morphism $g$ is bijective
on the set of closed points, it must be an
isomorphism. This completes the proof of the lemma.
\end{proof}

Note that as a consequence of Lemma \ref{basechange}, and the
fact that the morphism
$f$ is flat, the following holds: For any vector
bundle $E$ over $X$, the base change morphism
\begin{equation}\label{va.bci}
f^*F_{X*}E \,\longrightarrow\, F_{Z*}f^*E
\end{equation}
is an isomorphism, where $f$ is as in Eqn. \eqref{va.f}.

The following corollary is proved using the fact that the
homomorphism in Eqn. \eqref{va.bci} is an isomorphism.

\begin{corollary}\label{va.co.}
For any $f$ as in Eqn. \eqref{va.f} and
any integer $n\, \geq\, 1$, there is a canonical isomorphism
$$
f^*(F^n_X)_*E\,\longrightarrow\, (F^n_Z)_*f^*E\, . 
$$
\end{corollary}

\begin{proof}
We will use induction on $n$. For $n\, =\,1$, this is
the isomorphism in Eqn. \eqref{va.bci}. Assume that
we have constructed the isomorphism
\begin{equation}\label{va.bc2}
f^*(F^{n-1}_X)_*E \,\longrightarrow\, (F^{n-1}_Z)_*f^*E\, .
\end{equation}

Taking direct image by $F_Z$, the isomorphism in Eqn.
\eqref{va.bc2} gives an isomorphism
$$
F_{Z*}f^*(F^{n-1}_X)_*E\,\longrightarrow\, (F^{n}_Z)_*f^*E\, .
$$
On the other hand,
substituting $E$ by $(F^{n-1}_X)_*E$ in Eqn. \eqref{va.bci}
we get an isomorphism
$$
f^* (F^{n}_X)_*E\, =\, f^*F_{X*}(F^{n-1}_X)_*E
\,\longrightarrow\, F_{Z*}f^*(F^{n-1}_X)_*E\, .
$$
Combining the above two isomorphisms, an isomorphism
$$
f^*(F^n_X)_*E\,\longrightarrow\, (F^n_Z)_*f^*E
$$
is obtained. This completes the proof of the corollary.
\end{proof}

For any integer $n\, \geq\, 1$, there is a natural
homomorphism of ${\mathcal O}_{X}$ into the direct image
$(F^n_X)_*{\mathcal O}_X$. Let $B_n$ be the quotient bundle. So
we have a short exact sequence
\begin{equation}\label{eq55}
0\,\longrightarrow\, {\mathcal O}_{X}\,\longrightarrow\,
(F^n_X)_*{\mathcal O}_X\, \longrightarrow\, B_n
\,\longrightarrow\, 0
\end{equation}
of vector bundles over $X$. We note that the pullback of
this exact sequence by $F^n_X$ splits using the adjunction
homomorphism $(F^n_X)^*(F^n_X)_*{\mathcal O}_X\,\longrightarrow
\, {\mathcal O}_X$. Consequently,
\begin{equation}\label{eq551}
(F^n_X)^*(F^n_X)_*{\mathcal O}_X \, =\, {\mathcal O}_{X}\oplus
(F^n_X)^*B_n\, .
\end{equation}

Let $(B_n)_Z$ denote the quotient bundle
$(F^n_{Z})_*{\mathcal O}_Z/{\mathcal O}_{Z}$ over $Z$ which is
constructed by replacing $X$ by $Z$ in Eqn. \eqref{eq55}.

Corollary \ref{va.co.} gives the following:

\begin{corollary}\label{va.co2.}
For any $f$ as in Eqn. \eqref{va.f}, the
vector bundle $(B_n)_{Z}$ over $Z$ is canonically identified
with $f^* B_n$, where $n$ is any positive integer.
\end{corollary}

\begin{proof}
Set $E\, =\, {\mathcal O}_X$ in Corollary \ref{va.co.}.
The isomorphism in Corollary \ref{va.co.} evidently takes the line
subbundle
$$
{\mathcal O}_Z\, =\, f^*{\mathcal O}_X\, \subset\,
f^*(F^n_X)_*{\mathcal O}_X
$$
to the line subbundle
${\mathcal O}_Z\, \subset\, (F^n_Z)_*{\mathcal O}_Z$.
Hence the isomorphism in Corollary \ref{va.co.} for
$E\, =\, {\mathcal O}_X$ induces an isomorphism of $f^* B_n$
with $(B_n)_{Z}$. This completes the proof of the corollary.
\end{proof}

Let
$$
F^*_Z f^*F_{X*}E \,\longrightarrow\, F^*_Z F_{Z*}f^*E
$$
be the pull back of the isomorphism in Eqn. \eqref{va.bci}
by the morphism $F_Z$. On the other hand,
since $F_X\circ f\, =\, f\circ F_Z$, the vector bundle
$F^*_Zf^* F_{X*} E$ is identified with
$f^* F^*_X F_{X*} E$. Combining these
two isomorphisms we have an isomorphism
\begin{equation}\label{va.bc3}
f^* F^*_X F_{X*}E\,\longrightarrow\, F^*_Z F_{Z*} f^*E \, .
\end{equation}

The following lemma shows that the isomorphism in
Eqn. \eqref{va.bc3} is compatible with the filtration
of subbundles constructed in Eqn. \eqref{va.e4}.

\begin{lemma}\label{canfiletale}
Let $\{E_i\}_{i\geq 0}$ (respectively, $\{E'_i\}_{i\geq 0}$) be
the filtration of subbundles of the vector bundle $F^*_X F_{X*}E$
(respectively, $F^*_ZF_{Z*}f^*E$) constructed as in
Eqn. \eqref{va.e4}. For each $i\,\geq\, 0$, the isomorphism
in Eqn. \eqref{va.bc3} takes the subbundle
$f^*E_i\, \subset\, f^* F^*_X F_{X*}E$ to the subbundle $E'_i
\, \subset\, F^*_ZF_{Z*}f^*E$.
\end{lemma}

\begin{proof}
The following diagram is commutative
$$
\begin{matrix}
f^*F^*_XF_{X*}E & \stackrel{=}{\longrightarrow} &
F^*_ZF_{Z*}f^*E\\
\Big\downarrow &&\Big\downarrow \\
f^*E &= & f^*E
\end{matrix}
$$
where the vertical homomorphisms are defined as in
Eqn. \eqref{va.m1} and the top horizontal isomorphism
is constructed in Eqn. \eqref{va.bc3}. Consequently,
the isomorphism in Eqn. \eqref{va.bc3}
takes the subbundle $f^*E_1\, \subset\, f^*F^*_XF_{X*}E$ to
$E'_1\, \subset\, F^*_ZF_{Z*}f^*E$.

The Cartier connection on $F^*_XF_{X*}E$ induces a
connection on $f^*F^*_XF_{X*}E$. This induced connection
is taken to the Cartier connection on $F^*_ZF_{Z*}f^*E$
by the isomorphism in Eqn. \eqref{va.bc3}.

Since the filtration in Eqn. \eqref{va.e4} is
constructed from the subbundle $E_1$ and the Cartier connection
on $F^*_X F_{X*}E$, we
conclude that for each $i\,\geq\, 0$, the isomorphism in
Eqn. \eqref{va.bc3} takes the subbundle $f^*E_i$ to $E'_i$.
This completes the proof of the lemma.
\end{proof}

\begin{lemma}\label{etalevanish}
Let $X$ be an irreducible smooth projective
variety of dimension at least
three. Fix an ample line bundle $L$ over $X$.
There exist integers $r_1$ and $r_2$ such that
for each connected \'etale Galois covering
$f\,:\, Z\, \longrightarrow\, X$,
$$
H^i(Z,\, f^*L^{-r})\,=\,0
$$
provided $r\, > \, r_i$, where $i\, =\, 1,2$.
Furthermore, $r_1$ can be chosen to be
$M/{\rm degree}(L)$, where $M$ is given in
Proposition \ref{va.}.
\end{lemma}

\begin{proof}
Take any connected \'etale Galois covering $f\,:\,
Z\, \longrightarrow\, X$. Since $f$ is a finite morphism,
$$
H^{i}(Z,\, f^*L^{-r})\, =\, H^{i}(X,\, f_* f^*L^{-r})
$$
for all $i\, \geq\, 0$.
Combining this with the projection formula
$$
f_*f^*L^{-r} \,=\, f_*({\mathcal O}_{Z})\otimes L^{-r}
$$
we have
\begin{equation}\label{va.e.c}
H^i(Z,\, f^*L^{-r})\, =\,
H^i(X,\, f_*({\mathcal O}_{Z})\otimes L^{-r})
\end{equation}
for all $i\, \geq\, 0$.

For notational convenience, the vector bundle
$f_*{\mathcal O}_{Z}$ over $X$ will be denoted by $W$.

Tensoring the exact sequence in Eqn. \eqref{eq55} with
$W\otimes L^{-r}$ we have the short exact sequence
$$
0\,\longrightarrow\, W\otimes L^{-r}\,\longrightarrow\,
((F^n_X)_*{\mathcal O}_X) \otimes W\otimes L^{-r}\,\longrightarrow
\, W\otimes B_n\otimes L^{-r} \,\longrightarrow \,0\, .
$$
The long exact sequence of cohomologies for
this exact sequence gives the following:

If
\begin{equation}\label{vach1}
H^{i-1}(X,\, W\otimes B_n\otimes L^{-r})\, =\, 0
\end{equation}
and
\begin{equation}\label{vach2}
H^i(X,\, ((F^n_X)_*{\mathcal O}_X)\otimes W\otimes L^{-r})\, =\, 0\, ,
\end{equation}
then $H^i(X,\, W\otimes L^{-r}) \,=\, 0$.

The lemma will be proved by showing that for suitable values
of $r$, Eqn. \eqref{vach1} holds
for all $n$ and Eqn. \eqref{vach2} holds for all
sufficiently large $n$.

To prove that Eqn. \eqref{vach2} holds for $i\,=\,1,2$, first
note by the projection formula
$$
H^i(X,\, ((F^n_X)_*{\mathcal O}_X)\otimes W\otimes L^{-r})\, =\,
H^i(X,\, (F^n_X)_* (F^n_X)^* (W\otimes L^{-r}))
$$
for all $i\, \geq\, 0$. Now, since $F^n_X$ is a finite
morphism, we have
$$
H^i(X,\, (F^n_X)_* (F^n_X)^* (W\otimes L^{-r}))\, =\,
H^i(X,\, (F^n_X)^* (W\otimes L^{-r}))
$$
for all $i\, \geq\, 0$. Therefore, we have
\begin{equation}\label{vach5}
H^i(X,\, ((F^n_X)_*{\mathcal O}_X)\otimes W\otimes L^{-r})\, =\,
H^i(X,\, (F^n_X)^* (W\otimes L^{-r}))
\end{equation}
for all $i\, \geq\, 0$.

As the next step,
we will show that there is a natural isomorphism
$(F^n_X)^* W\,=\, W$ for all $n\, \geq\, 1$. To construct
this isomorphism, first note that since
the vector bundle $f^*W$ over $Z$ is trivializable,
there is a natural isomorphism
\begin{equation}\label{va.e.c2}
f^*W \, \longrightarrow\, F^*_Z f^*W\, .
\end{equation}
On the other hand, since
$F_X\circ f \, =\, f\circ F_Z$, there is a natural
isomorphism of $F^*_Z f^*W$ with $f^*F^*_X W$. This isomorphism
and the isomorphism in Eqn. \eqref{va.e.c2} together give
an isomorphism
\begin{equation}\label{156}
f^*W \, \longrightarrow\, f^*F^*_X W
\end{equation}
for all $n\, \geq\, 1$.

The isomorphism in Eqn. \eqref{156}
intertwines the actions of the Galois group
$\text{Gal}(f)$ on $f^*W$ and $f^*F^*_X W$.
Therefore, the
isomorphism in Eqn. \eqref{156} descends to an isomorphism
of $F^*_X W$ with $W$.
Since $F^*_X W\,=\, W$, we have
\begin{equation}\label{11f0}
(F^n_X)^* W\,=\, W\, .
\end{equation}

The isomorphism in Eqn. \eqref{vach5} and the
isomorphism in Eqn. \eqref{11f0} together give
\begin{equation}\label{va.e.cc}
H^i(X,\, ((F^n_X)_*{\mathcal O}_X)\otimes W\otimes L^{-r})\, =\,
H^i(X,\, W\otimes (F^n_X)^* L^{-r})
\end{equation}
(note that $(F^n_X)^* (W\otimes L^{-r})\, =\, ((F^n_X)^* W)\otimes
(F^n_X)^* L^{-r}\, =\, W\otimes (F^n_X)^* L^{-r}$). Since
$$
(F^n_X)^* L^{-r}\, =\, L^{-p^nr}\, ,
$$
where $p$ is the characteristic of the field $k$,
and $\dim X \, \geq\, 3$, using the
Serre vanishing it follows that
$$
H^i(X,\, W\otimes (F^n_X)^* L^{-r})\, =\,
H^i(X,\, W\otimes L^{-p^nr})\, =\, 0
$$
if $i\, =\, 1,2$ and $n$ is sufficiently
large \cite[page 111, Proposition 6.2.1]{EGA3}. Therefore, using Eqn.
\eqref{va.e.cc} we conclude that Eqn. \eqref{vach2} holds provided
$i\, =\, 1,2$ and $n$ is sufficiently large.

We will now prove that Eqn. \eqref{vach1} holds for $i\, =\, 1,2$
under appropriate conditions.

{}From the short exact sequence in Eqn. \eqref{eq551} it
follows that $L_{\rm max}((F^n_X)_*{\mathcal O}_X)\, \geq\, 0$ and
\begin{equation}\label{vach2p}
L_{\rm max}(B_n)\, \leq\, L_{\rm max}((F^n_X)_*{\mathcal O}_X)\, .
\end{equation}
We will use the ample line bundle $L$ to define degree
of a coherent sheaf on $X$.

Take any integer $r$ such that
$r\cdot \text{degree}(L)\, >\, M$, where $M$ is given in
Proposition \ref{va.}. From Proposition \ref{va.}
it follows that 
\begin{equation}\label{vach3}
L_{\rm max}((F^n_X)_*{\mathcal O}_X)\, <\,
r\cdot \text{degree}(L) \,=\, \text{degree}(L^{r})\, .
\end{equation}
Combining Eqn. \eqref{vach2p} and Eqn. \eqref{vach3} we conclude that
\begin{equation}\label{vach4}
L_{\rm max}(B_n\otimes L^{-r})\, =\,
L_{\rm max}(B_n)- \text{degree}(L^{r})\, \leq\,
L_{\rm max}((F^n_X)_*{\mathcal O}_X) - \text{degree}(L^{r})
\, <\, 0\, .
\end{equation}

Since $f$ is an \'etale Galois covering, the vector bundle
$f^*W$ over $Z$ is trivializable. Therefore, from
Eqn. \eqref{vach4} it follows that
$$
L_{\rm max}(f^*(W\otimes B_n\otimes L^{-r})) \, =\,
L_{\rm max}(f^*(B_n\otimes L^{-r})) \, <\, 0
$$
with respect to the polarization $f^*L$ on $Z$.
This inequality implies that
$$
H^0(Z,\, f^*(W\otimes B_n\otimes L^{-r}))\, =\, 0\, .
$$
Hence Eqn. \eqref{vach1} holds for all $n$ if $i\, =\, 1$ and
$r\, > \, M/\text{degree}(L)$.

Since both Eqn. \eqref{vach1} and Eqn. \eqref{vach2} hold for
$i\, =\,1$ and $r\, >\, M/\text{degree}(L)$, we know that
\begin{equation}\label{11n}
H^1(X, W\otimes L^{-r})\, =\, 0
\end{equation}
provided $r\, > \, M/\text{degree}(L)$. In view of Eqn.
\eqref{va.e.c}, this implies that
$$
H^1(Z,\, f^*L^{-r})\,=\,0
$$
if $r\, >\, r_1\, :=\, M/\text{degree}(L)$.

To prove Eqn. \eqref{vach1} for $i\, =\, 2$, consider
the exact sequence of sheaves
$$
0\,\longrightarrow\, F_{X*} {\mathcal O}_{X}\,\longrightarrow\,
F^{n}_{X*}{\mathcal O}_X\, \longrightarrow\, F_{X*} B_{n-1}
\,\longrightarrow\, 0
$$
obtained by taking the direct image of
Eqn. \eqref{eq55} by the morphism $F_X$ after substituting
$n-1$ for $n$ in Eqn. \eqref{eq55}. It gives an exact sequence
\begin{equation}\label{va3p1}
0\,\longrightarrow\, B_1\,\longrightarrow\,
B_n \, \stackrel{h_0}{\longrightarrow}\, F_{X*} B_{n-1}
\,\longrightarrow\, 0
\end{equation}
of vector bundles over $X$. Replacing $n$ by $n-1$ in Eqn.
\eqref{va3p1}, and then taking direct image by $F_X$, we have
$$
0\,\longrightarrow\, F_{X*} B_1 \,\longrightarrow\,
F_{X*}B_{n-1} \, \stackrel{h_1}{\longrightarrow}\,
F_{X*} F_{X*} B_{n-2} \,\longrightarrow\, 0\, .
$$
More generally, replacing $n$ by $n-j$ in Eqn. \eqref{va3p1},
and then taking direct image by $F^j_X$, we have
$$
0\,\longrightarrow\, F^j_{X*} B_1 \,\longrightarrow\,
F^j_{X*}B_{n-j} \, \stackrel{h_j}{\longrightarrow}
\, F^{j+1}_{X*}B_{n-j-1}\,\longrightarrow\, 0\, ,
$$
where $0\, \leq\, j\,\leq\, n-1$.

For $j\, \in\, [0\, ,n-1]$, set $U_j\, :=\, F^j_{X*}B_{n-j}$.
Also, set $U_n\, =\, 0$.
We have a filtration of quotients
\begin{equation}\label{11fil}
B_n\, :=\, U_0 \, \stackrel{h_0}{\longrightarrow}\, U_1\,
\stackrel{h_1}{\longrightarrow}\,
U_2\, \stackrel{h_2}{\longrightarrow}\,\cdots\,
\stackrel{h_{n-2}}{\longrightarrow}\, U_{n-1}\,
\stackrel{h_{n-1}}{\longrightarrow}
\, U_n\, =\, 0
\end{equation}
of the vector bundle $B_n$, where the homomorphisms $h_j$,
$j\, \in\, [0\, ,n-1]$, are constructed above. Note that
$\text{kernel}(h_j) \, =\, F^j_{X*} B_1$, $j\, \in\, [0\, , n-1]$,
with the convention that $(F^0_{X})_* B_1 \,=\, B_1$ (recall
that $F^0_X$ denotes the identity morphism of $X$). 

Using the filtration in Eqn. \eqref{11fil} we conclude that
Eqn. \eqref{vach1} holds for $i\, =\, 2$ if
\begin{equation}\label{11f1}
H^1(X,\, W\otimes L^{-r} \otimes (F^j_X)_* B_1)\, =\, 0
\end{equation}
for all $j\, \in\, [0\, , n-1]$.

The projection formula says that
$W\otimes L^{-r}\otimes (F^j_X)_* B_1\, =\,
(F^j_X)_*((F^j_X)^*(W\otimes L^{-r})\otimes B_1)$.
We have $(F^j_X)^* W\, =\,W$
(see Eqn. \eqref{11f0}) and $(F^j_X)^* L^{-r}\, =\,
L^{-p^jr}$. Using these and the
fact that $F^j_X$ is a finite morphism we have
\begin{equation}\label{11f2}
H^1(X,\, W\otimes L^{-r}\otimes (F^j_X)_* B_1)\, =\,
H^1(X,\, (F^j_X)_*((F^j_X)^*(W\otimes L^{-r})\otimes B_1))
\end{equation}
$$
=\, H^1(X,\, (F^j_X)^*(W\otimes L^{-r})\otimes B_1)\, =\,
H^1(X,\, W\otimes L^{-rp^j}\otimes B_1)\, .
$$

There is a positive integer $m$, such that there
exists a surjective homomorphism
$$
\psi\, :\, (L^{-m})^{\oplus m} \, \longrightarrow\, B^*_1\, .
$$
Given such a homomorphism $\psi$, consider the exact
sequence of vector bundles
\begin{equation}\label{2}
0\, \longrightarrow\, B_1\, \stackrel{\psi^*}{\longrightarrow}
\, (L^{m})^{\oplus m} \, \longrightarrow\,
{\mathcal H} \, \longrightarrow\, 0
\end{equation}
over $X$.

Fix an integer $r_2$ satisfying the two conditions:
\begin{enumerate}
\item{} $r_2\cdot\text{degree}(L)
\, \geq\, L_{\rm max}(\mathcal H)$, and

\item{} $r_2 \, \geq\, m+ r_1$,
where $r_1$ is defined in the statement of the lemma.
\end{enumerate}

We will show that Eqn. \eqref{vach1} with $i\, =\, 2$
holds for all $n$ and all $r\, >\, r_2$. For that purpose
we will show that
\begin{equation}\label{01}
H^1(X,\, W \otimes B_1\otimes L^{-r}) \, =\, 0
\end{equation}
for all $r\, >\, r_2$.

Take any $r\, >\, r_2$, where $r_2$ is the above
integer. Let
\begin{equation}\label{11a}
\longrightarrow\,
H^0(X,\, W \otimes{\mathcal H}\otimes L^{-r})\, \longrightarrow\,
H^1(X,\, W \otimes B_1\otimes L^{-r}) \, \longrightarrow\,
H^1(X,\, W \otimes L^{m-r})^{\oplus m}
\end{equation}
be the long exact sequence of cohomologies for the
short exact sequence of sheaves obtained by tensoring the exact
sequence in Eqn. \eqref{2} with $W\otimes L^{-r}$.

We already proved that $H^1(X, W\otimes L^{-r'})\, =\, 0$ if
$r'\, > \, r_1$; see Eqn. \eqref{11n}. Therefore,
\begin{equation}\label{4}
H^1(X,\, W\otimes L^{m-r})\,=\, 0
\end{equation}
as $r-m\, >\, r_2-m \, \geq\, r_1$
(see the definition of $r_2$ given above).

On the other hand, as $L_{\text{max}}(W)\, =\, 0$
(the pull back $f^*W$ is trivializable), we have
$$
L_{\text{max}}(W\otimes{\mathcal H}\otimes L^{-r})\, =\,
L_{\text{max}}({\mathcal H}\otimes L^{-r})
\, =\, L_{\text{max}}({\mathcal H}) -\text{degree}(L^r)
\, <\, 0
$$
(the last inequality follows from the conditions that
$r\, >\, r_2$ and $r_2\cdot\text{degree}(L)
\, \geq\, L_{\rm max}(\mathcal H)$). Therefore, we have
\begin{equation}\label{3}
H^0(X,\, W \otimes {\mathcal H}\otimes L^{-r})\, =\, 0\, .
\end{equation}

Combining Eqn. \eqref{4} and Eqn. \eqref{3}
with Eqn. \eqref{11a} we conclude that Eqn. \eqref{01}
holds. In view of Eqn. \eqref{11f1} and Eqn. \eqref{11f2},
this implies that Eqn. \eqref{vach1} with $i\, =\, 2$ holds
for all $n$. We showed earlier that Eqn. \eqref{vach2}
$i\, =\, 2$ holds for all $n$ sufficiently large.

Therefore, if $r\, \geq\, r_2$, then $H^2(X,\, W\otimes L^{-r})
\,=\, 0$. In view of Eqn.
\eqref{va.e.c}, this completes the proof of the lemma.
\end{proof}

We will need the following lemma for our purpose.

\begin{lemma}\label{vanish3}
Let $X$ and $L$ be as in Lemma \ref{etalevanish}. Let $V$ be a
vector bundle over $X$. Then there is an integer $N(V)$
such for any $r\,> \, N(V)$, and for each connected \'etale Galois
cover $f\,:\, Z\, \longrightarrow\, X$,
$$
H^1(Z,\, (F^n_Z)^*f^*(V\otimes L^{-r}))\,=\, 0
$$
for all $n\, \geq\, 0$.
\end{lemma}

\begin{proof}
There is an integer $m$ such that there is a
surjective homomorphism
$$
\psi\, :\, (L^{-m})^{\oplus m} \, \longrightarrow\, V^*
$$
of vector bundles. The homomorphism $\psi$
gives an exact sequence of vector bundles
\begin{equation}\label{11ex}
0\, \longrightarrow\, V\, \stackrel{\psi^*}{\longrightarrow}
\, (L^{m})^{\oplus m}\, \longrightarrow\,
{\mathcal A} \, \longrightarrow\, 0
\end{equation}
over $X$.

Fix an integer $N(V)$ satisfying the two conditions:
\begin{enumerate}
\item{} $N(V)\cdot\text{degree}(L)
\, \geq \, L_{\rm max}(\mathcal A)$, and

\item{} $N(V) \, \geq\, m+ r_1$,
where $r_1$ is prescribed in Lemma \ref{etalevanish}.
\end{enumerate}
These conditions ensure that if $r\, > \, N(V)$, then
\begin{equation}\label{160}
L_{\rm max}({\mathcal A}\otimes L^{-r})\,<\, 0
\end{equation}
and
\begin{equation}\label{161}
H^1(Z,\, f^*L^{p^n(m-r)})\, =\, 0
\end{equation}
for all $n\,\geq\, 0$ (see Lemma \ref{etalevanish}). Since
$$
(F^n_Z)^*f^* L^{m-r} \, =\, f^*L^{p^n(m-r)}\, ,
$$
{}from Eqn. \eqref{161} we conclude that
\begin{equation}\label{161p}
H^1(Z,\, (F^n_Z)^*f^* L^{m-r})\, =\, 0
\end{equation}
for all $n\,\geq\, 0$.

Take any $r\, > \, N(V)$. Take any $f$ as in the statement of the
lemma. From Eqn. \eqref{160} we have
$$
L_{\rm max}(f^*({\mathcal A}\otimes L^{-r}))\,<\, 0
$$
with respect to the polarization on $Z$ obtained by pulling back
the polarization on $X$. This implies that
$$
L_{\rm max}((F^n_Z)^* f^*({\mathcal A}\otimes L^{-r}))
\, =\, p^n L_{\rm max}(f^*({\mathcal A}\otimes L^{-r}))\,<\, 0
$$
for all $n\, \geq\, 0$. Consequently,
\begin{equation}\label{112}
H^0(Z, \, (F^n_Z)^* f^*({\mathcal A}\otimes L^{-r}))\, =\, 0
\end{equation}
and for all $n\, \geq\, 0$.

Tensor the exact sequence in Eqn. \eqref{11ex} with
$L^{-r}$ and then pull it back by the morphism $f\circ F^n_Z$.
This gives the exact sequence of vector bundles
$$
0\, \longrightarrow\, (F^n_Z)^* f^*(V\otimes L^{-r})\,\longrightarrow
\, ((F^n_Z)^*f^* L^{m-r})^{\oplus m}\, \longrightarrow\,
(F^n_Z)^*f^* ({\mathcal A}\otimes L^{-r}) \, \longrightarrow\, 0
$$
over $Z$.
Consider the corresponding long exact sequence of cohomologies
$$
H^0(Z,\, (F^n_Z)^*f^* ({\mathcal A}\otimes L^{-r})) \,
\longrightarrow \, H^1(Z,\, (F^n_Z)^* f^*(V\otimes L^{-r})) \,
\longrightarrow \,
H^1(Z,\, (F^n_Z)^*f^* L^{m-r})^{\oplus m}\, .
$$
In view of Eqn. \eqref{112} and Eqn. \eqref{161p} it implies that
$$
H^1(Z,\, (F^n_Z)^* f^*(V\otimes L^{-r}))\, =\, 0\, .
$$
This completes the proof of the lemma.
\end{proof}

The following lemma is very similar to Lemma \ref{vanish3}.

\begin{lemma}\label{vanish3-1}
Let $X$ and $L$ be as in Lemma \ref{etalevanish}. Fix vector
bundles $V_1, V_2, \cdots , V_\ell$ over $X$.
Then there is an integer $N\, =\, N(V_1,\cdots , V_\ell)$
such that
$$
H^1(Z, \, (\bigotimes_{j=0}^{n-1}(F^{j}_Z)^*f^* V_{n_j})
\otimes f^* L^{-p^nr}) \, =\, 0
$$
provided $r\, >\, N$,
where $n$ is any positive integer and $n_j\, \in\,
[1\, ,\ell]$ for all $j \, \in\, [0\, ,n-1]$, and
$f\, :\, Z\, \longrightarrow\, X$ is any connected
\'etale Galois cover.
\end{lemma}

\begin{proof}
As in Eqn. \eqref{11ex}, there is an integer $m$ and exact
sequences
\begin{equation}\label{162}
0\, \longrightarrow\, V_i\, \stackrel{\psi_i}{\longrightarrow}
\,(L^m)^{\oplus m} \, \longrightarrow\,
{\mathcal A}_i \, \longrightarrow\, 0
\end{equation}
for all $i\, \in\, [1\, , \ell]$.
For any positive integer $n$ and integers $n_j\, \in\, [1\, ,\ell]$,
where $j\, \in\, [0\, ,n-1]$, the homomorphisms $\psi_i$
in Eqn. \eqref{162} together give an injective homomorphism
of vector bundles
\begin{equation}\label{165}
\bigotimes_{j=0}^{n-1} (F^{j}_Z)^*f^* V_{n_j}
\, \longrightarrow\, (f^*L^{m(p^n-1)/(p-1)})^{\oplus m^n}\, .
\end{equation}
The quotient bundle for
the homomorphism in Eqn. \eqref{165} is filtered by subbundles
such that each successive quotient is of the form
$$
(\bigotimes_{j\in I} (F^{j}_Z)^*f^* V_{n_j})\otimes
(\bigotimes_{j\in I^c}(F^{j}_Z)^*f^*{\mathcal A}_{n_j})\, ,
$$
where $I$ is some subset of $[0\, ,n-1]$. Therefore, using
this filtration, from the
long exact sequence of cohomologies for the short exact sequence
of vector bundles obtain from Eqn. \eqref{165} tensored with
$L^{-p^nr}$ we conclude that if
\begin{equation}\label{166}
H^0(Z,\, (\bigotimes_{j\in I} (F^{j}_Z)^*f^* V_{n_j})\otimes
(\bigotimes_{j\in I^c}(F^{j}_Z)^*f^*{\mathcal A}_{n_j})
\otimes f^* L^{-p^nr})\, =\, 0
\end{equation}
for all $I\, \subset\, [0\, ,n-1]$, and
\begin{equation}\label{167}
H^1(Z,\, f^*L^{m(p^n-1)/(p-1)}\otimes f^* L^{-p^nr})\, =\, 0\, ,
\end{equation}
then
$$
H^1(Z, \, (\bigotimes_{j=0}^{n-1}(F^{j}_Z)^*f^* V_{n_j})
\otimes f^* L^{-p^nr})\, =\, 0\, .
$$

Take a nonnegative integer $N\, :=\, N(V_1, \cdots , V_\ell)$
satisfying the following conditions:
\begin{enumerate}
\item{} $N\cdot\text{degree}(L) \, \geq\,
L_{\rm max}({\mathcal A}_i)$ and $N\cdot\text{degree}(L)
\, \geq\, L_{\rm max}(V_i)$ for all $i\, \in\,
[1\, , \ell]$;
\item{} $N\, \geq\, r_1+ m$, where $r_1$ is defined
in Lemma \ref{etalevanish} and $m$ is as in Eqn. \eqref{162}.
\end{enumerate}

If $r\, \,>\, N$, then from the given condition
$N\, \geq\, r_1+ m$ it follows that
$$
p^nr - m(p^n-1)/(p-1) \, >\, r_1\, .
$$
Therefore, Eqn. \eqref{167} follows from Lemma \ref{etalevanish}
provided $r\, \,>\, N$.

To prove Eqn. \eqref{166}, we note that
$$
L_{\rm max}((\bigotimes_{j\in I} (F^{j}_Z)^*f^* V_{n_j})\otimes
(\bigotimes_{j\in I^c}(F^{j}_Z)^*f^*{\mathcal A}_{n_j})
\otimes f^* L^{-p^nr})
\,=\,(\sum_{j\in I}L_{\rm max}((F^{j}_Z)^*f^* (V_{n_j}\otimes L^{-r}))
$$
$$
+ \sum_{j\in I^c }L_{\rm max}((F^{j}_Z)^*f^*({\mathcal A}_{n_j}
\otimes L^{-r})) + \frac{r(p^n-1)}{p-1} -rp^n)\text{degree}(f)
\, \leq \,(\frac{r(p^n-1)}{p-1} -rp^n)\text{degree}(f)
$$
as $L_{\rm max}(V_{i}\otimes L^{-r}), L_{\rm max}({\mathcal A}_i
\otimes L^{-r})\, \leq\, 0$ for all $i\, \in\, [1\, ,\ell]$
(this is the first of the two conditions that define $N$).
Therefore, as $r(p^n-1)/(p-1) \, <\, rp^n$, we conclude that
Eqn. \eqref{166} holds if $r\, \,>\, N$.
This completes the proof of the lemma.
\end{proof}

\begin{proposition}\label{vanish-imp}
Let $X$ and $L$ be as in Lemma \ref{etalevanish}.
There is an integer $N$ such that for all
$r\,>\, N$ the following holds:
$$
H^1(Z, (F^n_Z)^*(F^n_Z)_* f^*L^{-p^nr})\,=\, 0
$$
for all $n\, \geq\, 0$ and all connected \'etale Galois cover
$f\,:\, Z\, \longrightarrow\, X$.
\end{proposition}

\begin{proof}
Set $E\, :=\,(F^{n-1}_Z)_* {\mathcal O}_Z$. Consider
the filtration of subbundles of
$F^*_Z F_{Z*} E$ constructed as in Eqn.
\eqref{va.e4}. Pulling this filtration back by $F^{n-1}_Z$
we get a filtration of subbundles of
$$
(F^{n-1}_Z)^* F^*_Z F_{Z*} E
\, =\,(F^n_Z)^*(F^n_Z)_* {\mathcal O}_Z
$$
whose each successive quotient is the subbundle
$$
(F^{n-1}_Z)^* (E\otimes {\mathcal G}_i)\, \subset\,
(F^{n-1}_Z)^* (E\otimes \Omega^{\otimes i}_X)\, ,
$$
where $0 \,\leq\,i\, \leq\, (p-1){\dim}(X)$ and
${\mathcal G}_i$ is the vector bundle on $Z$
defined in Proposition \ref{canfil}(3).

In view of Lemma \ref{canfiletale}
and the fact that the connection on
$f^*F^*_XF_{X*}{\mathcal O}_X$
induced by the Cartier connection on $F^*_XF_{X*}{\mathcal O}_X$
is taken to the Cartier connection on $F^*_ZF_{Z*}{\mathcal O}_Z$
by the isomorphism in Eqn. \eqref{va.bc3}
(see the proof of Lemma \ref{canfiletale}), we conclude that
${\mathcal G}_i\,
=\, f^*{\mathcal G}^X_i$, where ${\mathcal G}^X_i\,
\subset\, \Omega^{\otimes i}_X$ is the subbundle
defined in Proposition \ref{canfil}(3). Therefore,
\begin{equation}\label{168}
(F^{n-1}_Z)^* (E\otimes {\mathcal G}_i)\, =\,
((F^{n-1}_Z)^* E)\otimes (F^{n-1}_Z)^* f^*{\mathcal G}^X_i
\, =\, ((F^{n-1}_Z)^* (F^{n-1}_Z)_* {\mathcal O}_Z)
\otimes (F^{n-1}_Z)^* f^*{\mathcal G}^X_i\, .
\end{equation}

Repeating the above argument after replacing $n$ by $n-1$,
for all $i\, \in\, [0\, ,(p-1){\dim}(X)]$ we get a
filtration of subbundles of
the vector bundle in Eqn. \eqref{168}
whose successive quotients are of the form
$$
((F^{n-2}_Z)^* (F^{n-2}_Z)_* {\mathcal O}_Z)
\otimes (F^{n-2}_Z)^* f^*{\mathcal G}^X_j
\otimes (F^{n-1}_Z)^* f^*{\mathcal G}^X_i\, ,
$$
where $0 \,\leq\, j\, \leq\, (p-1){\dim}(X)$.

In fact, iterating the above argument we conclude the
following:

To prove the proposition it is enough to show that
there is an integer $N$ such that for any
$r\,>\, N$,
\begin{equation}\label{169}
H^1(Z,\, (\bigotimes_{i=0}^{n-1}(F^{i}_Z)^*f^*
{\mathcal G}^X_{n_i})
\otimes f^* L^{-p^nr}) \, =\, 0 
\end{equation}
for all nonnegative integer $n$ and all $n_i\, \in\, [0\, ,
(p-1){\dim}(X)]$.

Setting $\ell\, =\, (p-1){\dim}(X)+1$ and
$V_j\, =\, {\mathcal G}^X_{j-1}$ in Lemma \ref{vanish3-1}
we get an integer $N$ such that Eqn. \eqref{169} holds for
all $r\,>\, N$. This completes the proof of the proposition.
\end{proof}

\begin{proposition}\label{vanish4}
Let $X$ be an irreducible smooth projective
variety of dimension $d$, with $d\, \geq\, 3$.
Fix an ample line bundle $L$ over $X$.
There is an integer $N\, =\, N(X,L)$ with the following
property. Let $f\, :\, Z\,\longrightarrow\, X$ be a
connected \'etale Galois cover and $E$ an $F$--trivial vector
bundle over $Z$. Then
\begin{enumerate}
\item{} $H^1(Z,\, E\otimes f^*L^{-r})\, =\, 0$ and

\item{} $H^1(Z,\, E\otimes \Omega_Z\otimes f^*L^{-r}) \, =\, 0$
\end{enumerate}
for all $r\, >\, N$.
\end{proposition}

\begin{proof}
Take $f$ and $E$ as in the statement of the proposition.
Let $n$ be a positive integer such that the vector bundle
$(F^n_Z)^* E$ is trivializable.
To prove statement (1), let
\begin{equation}\label{16a}
0\,\longrightarrow\, {\mathcal O}_{Z}\,\longrightarrow\,
(F^n_{Z})_*{\mathcal O}_Z\, \longrightarrow\, (B_n)_Z
\,\longrightarrow\, 0
\end{equation}
be the exact sequence of vector bundles over $Z$
constructed as in Eqn. \eqref{eq55}. Consider the long
exact sequence of cohomologies for the
short exact sequence obtained
by tensoring Eqn. \eqref{16a} with $E\otimes f^*L^{-r}$. From
it we conclude that if
\begin{equation}\label{16b}
H^0(Z,\, (B_n)_Z\otimes E\otimes f^*L^{-r}) \, =\, 0
\end{equation}
and
\begin{equation}\label{16c}
H^1(Z,\, E\otimes f^*L^{-r}\otimes (F^n_{Z})_*{\mathcal O}_Z)
\, =\, 0\, , 
\end{equation}
then $H^1(Z,\, E\otimes f^*L^{-r})\, =\, 0$.

Using the projection formula and the fact that the
Frobenius morphism $F_Z$ is flat we have
$$
H^1(Z,\, E\otimes f^*L^{-r}\otimes (F^n_{Z})_*{\mathcal O}_Z)
\,=\, H^1(Z,\, (F^n_{Z})_* (F^n_Z)^*(E\otimes f^*L^{-r}))
$$
$$
\, =\, H^1(Z,\, (F^n_Z)^*(E\otimes f^*L^{-r})) \,=\,
H^1(Z,\, ((F^n_Z)^*E)\otimes f^*L^{-p^nr})\, \cong\,
H^1(Z,\, f^*L^{-p^nr})^{\oplus r_E} \, ,
$$
where $r_E\, =\, \text{rank}(E)$; recall that the
vector bundle $(F^n_Z)^* E$ is trivializable.

Lemma \ref{etalevanish} says that
$H^1(Z,\,  f^*L^{-p^nr})\, =\, 0$
if $p^nr\, >\, r_1$. Therefore, Eqn. \eqref{16c}
holds if $r\, >\, r_1$, where $r_1$ is given by Lemma
\ref{etalevanish}.

To prove Eqn. \eqref{16b} we first recall Eqn. \eqref{vach2p}
which says that
$$
L_{\rm max}((B_n)_Z)\, \leq\,
L_{\rm max}((F^n_Z)_*{\mathcal O}_Z)\, .
$$
Since $L_{\rm max}(E)\, =\, 0$, this inequality gives
\begin{equation}\label{16d}
L_{\rm max}((B_n)_Z\otimes E\otimes f^*L^{-r}) =
L_{\rm max}((B_n)_Z) - \text{degree}(f^*L^{r})
{\leq} L_{\rm max}((F^n_Z)_*{\mathcal O}_Z) -
\text{degree}(f^*L^{r})\, .
\end{equation}

Corollary \ref{va.co.} implies that
$$
L_{\rm max}((F^n_Z)_*{\mathcal O}_Z)\, =\,
L_{\rm max}(f^*(F^n_X)_*{\mathcal O}_X)\, .
$$
Hence from Proposition \ref{va.} we conclude that the
right--hand side of Eqn. \eqref{16d} is negative if
$r\, >\, r_1$, where $r_1$ is given in Lemma
\ref{etalevanish}.

Consequently, Eqn. \eqref{16b} holds for all $r\, >\, r_1$.
This completes the proof of statement (1) of the
proposition.

Statement (2) is proved in a similar way.

Tensoring Eqn. \eqref{16a} with $E\otimes\Omega_Z\otimes
f^*L^{-r}$ and considering the corresponding long
exact sequence of cohomologies we conclude that if
\begin{equation}\label{16e}
H^0(Z,\, (B_n)_Z\otimes E\otimes\Omega_Z\otimes f^*L^{-r})
\, =\, 0
\end{equation}
and
\begin{equation}\label{16f}
H^1(Z,\, E\otimes \Omega_Z\otimes f^*L^{-r}\otimes
(F^n_{Z})_*{\mathcal O}_Z)\, =\, 0\, , 
\end{equation}
then $H^1(Z,\, E\otimes\Omega_Z\otimes f^*L^{-r})\, =\, 0$.

We have
$$
L_{\rm max}((B_n)_Z\otimes E\otimes\Omega_Z\otimes f^*L^{-r}) 
\,=\, L_{\rm max}((B_n)_Z\otimes E\otimes f^*L^{-r})
+L_{\rm max}(\Omega_Z)
$$
$$
=\, L_{\rm max}((B_n)_Z\otimes
E\otimes f^*L^{-r})+ L_{\rm max}(f^*\Omega_X)\, .
$$
We saw that $L_{\rm max}((B_n)_Z\otimes
E\otimes f^*L^{-r})\, <\, 0$ if $r\, >\,  r_1$.
Therefore,
$$
L_{\rm max}((B_n)_Z\otimes E\otimes\Omega_Z\otimes f^*L^{-r})
\, <\, 0
$$
for all $r\, >\, r_1+L_{\rm max}(\Omega_X)/\text{degree}(L)$.
Hence Eqn. \eqref{16e} holds if
$$
r\,>\, r_1+L_{\rm max}(\Omega_X)/\text{degree}(L)\, .
$$

As before, using the projection formula and the fact that the
Frobenius morphism $F_Z$ is flat we have
$$
H^1(Z,\, E\otimes \Omega_Z\otimes f^*L^{-r}\otimes
(F^n_{Z})_*{\mathcal O}_Z)
\, =\, H^1(Z,\, (F^n_Z)^*(E\otimes\Omega_Z\otimes f^*L^{-r}))
$$
$$
\cong\, H^1(Z,\, (F^n_Z)^*(\Omega_Z\otimes
f^*L^{-r}))^{\oplus r_E}\,=\, H^1(Z,\,
(F^n_Z)^*f^*(\Omega_X\otimes L^{-r}))^{\oplus r_E}\, .
$$

Now,
Lemma \ref{vanish3} says that there is an integer $N'$ such
that
$$
H^1(Z,\, (F^n_Z)^*f^*(\Omega_X\otimes L^{-r}))\, =\, 0
$$
for all $r\, >\, N'$. Therefore, Eqn. \eqref{16f} holds if
$r\, >\, N'$.

Thus, statement (2) of the proposition holds if $r\, >\,
\text{Max}\{N',r_1+ L_{\rm max}(\Omega_X/\text{degree}(L))\}$.
This completes the proof of the proposition.
\end{proof}

Proposition \ref{vanish4}(2) will be very useful in Section
\ref{sec.-inj.}.

\begin{lemma}\label{vanish6}
Let $X$ and $L$ be as in Proposition \ref{vanish4}.
There is an integer $N\, =\, N(X,L)$ with the following
property. Let $f\, :\, Z\,\longrightarrow\, X$ be any
connected \'etale Galois cover and $E$ an $F$--trivial vector
bundle over $Z$. Then
$$
H^2(Z,\, E\otimes f^*L^{-r}) \, =\, 0
$$
for all $r\,>\, N$.
\end{lemma}

\begin{proof}
Take any $f$ and $E$ as in the statement of the lemma.
Let $n$ be an integer such that the vector
bundle $(F^n_Z)^*E$ is trivializable.

Let
$$
H^1(Z,\, (B_n)_{Z}\otimes E\otimes f^*L^{-r})\,\longrightarrow
\, H^2(Z,\, E\otimes f^*L^{-r}) \,\longrightarrow\,
H^2(Z,\, E\otimes ((F^n_{Z})_*{\mathcal O}_Z)\otimes
f^* L^{-r})
$$
be the exact sequence of cohomologies for the short exact
sequence of vector bundles obtained by tensoring the
exact sequence in Eqn. \eqref{16a} with $E\otimes
f^* L^{-r}$. Using this exact sequence it follows that if
\begin{equation}\label{16g}
H^2(Z,\, E\otimes ((F^n_{Z})_*{\mathcal O}_Z)\otimes
f^* L^{-r})\, =\, 0
\end{equation}
and
\begin{equation}\label{16h}
H^1(Z,\, (B_n)_{Z}\otimes E\otimes f^*L^{-r})\, =\, 0\, ,
\end{equation}
then $H^2(Z,\, E\otimes f^*L^{-r})\, =\, 0$.

Using the projection formula and the fact that the
Frobenius morphism $F_Z$ is flat we have
$$
H^2(Z,\, E\otimes ((F^n_{Z})_*{\mathcal O}_Z)\otimes
f^* L^{-r})\, =\,
H^2 (Z,\, (F^n_Z)^*(E\otimes f^*L^{-r}))\, \cong\,
H^2(Z,\, f^*L^{-p^nr})^{\oplus r_E}\, ,
$$
where $r_E\, =\, \text{rank}(E)$. Now,
Lemma \ref{etalevanish} says that
$H^2(Z,\, f^*L^{-p^nr})\, =\, 0$
if $r\, >\, r_2$, where $r_2$ is given by
Lemma \ref{etalevanish}.

Therefore, Eqn. \eqref{16g} holds if $r\, >\, r_2$.

To prove Eqn. \eqref{16h}, tensor the exact sequence in Eqn.
\eqref{16a} with $(B_n)_{Z}\otimes E\otimes f^*L^{-r}$ and
consider the corresponding long exact sequence of
cohomologies. From it we conclude that Eqn. \eqref{16h} holds if
\begin{equation}\label{16i}
H^0(Z,\, (B_n)_Z\otimes (B_n)_Z\otimes E\otimes f^*L^{-r})
\, =\, 0
\end{equation}
and
\begin{equation}\label{16j}
H^1(Z,\,((F^n_{Z})_*{\mathcal O}_Z)\otimes (B_n)_Z\otimes E\otimes
f^*L^{-r})\, =\, 0\, .
\end{equation}

Note that
\begin{equation}\label{16k}
L_{\rm max}((B_n)_Z\otimes (B_n)_Z\otimes E\otimes f^*L^{-r})
\,=\, 2L_{\rm max}((B_n)_Z)- r\cdot \text{degree}(f^*L)
\end{equation}
as $L_{\rm max}(E)\, =\, 0$ and
$L_{\rm max}((B_n)_Z\otimes (B_n)_Z)\, =\, 2L_{\rm max}((B_n)_Z)$.

Corollary \ref{va.co2.} says that $(B_n)_Z\, =\, f^*B_n$.
Using Eqn. \eqref{vach2p} and Proposition \ref{va.} we know
that
$$
L_{\rm max}(B_n)\,\leq\, M\, ,
$$
where $M$ is given in Proposition \ref{va.}. In other words,
$$
2L_{\rm max}(B_n)- r\cdot \text{degree}(L)\, <\, 0
$$
if $r\, >\, 2M/\text{degree}(L)$. Since $(B_n)_Z\, =\, f^*B_n$,
this implies that
$$
2L_{\rm max}((B_n)_Z)- r\cdot\text{degree}(f^*L) \, <\, 0
$$
if $r\, >\, 2M/\text{degree}(L)$.

Therefore, using Eqn. \eqref{16k} we conclude that Eqn. \eqref{16i}
holds if $r\, >\, 2M/\text{degree}(L)$.

To prove Eqn. \eqref{16j}, first note that using projection
formula and the fact that
$(F^n_{Z})^*E$ is trivializable we have
$$
H^1(Z,\,((F^n_{Z})_*{\mathcal O}_Z)\otimes (B_n)_Z\otimes E\otimes
f^*L^{-r})\, =\, H^1(Z,\, (F^n_{Z})_*(F^n_{Z})^*(
(B_n)_Z\otimes E\otimes f^*L^{-r}))
$$
$$
=\, H^1(Z,\,(F^n_{Z})^*((B_n)_Z\otimes E\otimes
f^*L^{-r}))\, \cong\, H^1(Z,\,(F^n_{Z})^*((B_n)_Z\otimes
f^*L^{-r}))^{r_E}\, ,
$$
where $r_E\, =\, \text{rank}(E)$. Therefore, Eqn. \eqref{16j}
holds if
\begin{equation}\label{16m}
H^1(Z,\,(F^n_{Z})^*((B_n)_Z\otimes f^*L^{-r}))\, =\, 0\, .
\end{equation}

Tensor the exact sequence in Eqn. \eqref{16a} with
$f^*L^{-r}$ and then pull it back by $F^n_{Z}$. Consider
the long exact sequence of cohomologies
\begin{equation}\label{16n}
\longrightarrow\, H^1(Z,\,(F^n_{Z})^*(((F^n_{Z})_*{\mathcal O}_Z)
\otimes f^*L^{-r})) \, \longrightarrow\, H^1(Z,\,
(F^n_{Z})^*((B_n)_Z\otimes f^*L^{-r}))
\end{equation}
$$
\longrightarrow\, H^2(Z,\, (F^n_{Z})^* f^*L^{-r})
\,\longrightarrow
$$
corresponding to the resulting exact sequence of vector bundles.

Using projection formula we have
$$
(F^n_{Z})^*(((F^n_{Z})_*{\mathcal O}_Z) \otimes
f^*L^{-r})\,=\, (F^n_{Z})^*(F^n_{Z})_*(F^n_{Z})^*f^*L^{-r}
\,=\, (F^n_{Z})^*(F^n_{Z})_*f^*L^{-p^nr}\, .
$$
Therefore, from Proposition \ref{vanish-imp} we conclude
that there is an integer $N'$ such that
$$
H^1(Z,\,(F^n_{Z})^*(((F^n_{Z})_*{\mathcal O}_Z) \otimes
f^*L^{-r}))\,=\, H^1(Z,\, (F^n_{Z})^*(F^n_{Z})_*f^*L^{-p^nr})
\, =\, 0
$$
if $r\, >\, N'$. In view of this, using Lemma
\ref{etalevanish}, from Eqn. \eqref{16n} we conclude that
there is an integer $N''$ such that Eqn. \eqref{16m}
holds if $r\, >\, N''$. This completes the proof of
the proposition.
\end{proof}

\begin{proposition}\label{vi}
Let $X$ be an irreducible smooth projective
variety of dimension $d$, with $d\, \geq\, 3$.
Fix an ample line bundle $L$ over $X$. There is an integer
$N\, =\, N(X,L)$ with the following property.
Let $f\, :\, Z\,\longrightarrow\, X$ be a connected
\'etale Galois cover and $E$ an $F$--trivial vector
bundle over $Z$. Then for any smooth divisor
$D\, \in\, \vert L^r\vert$ with $r\, >\, N$,
$$
H^1(f^{-1}(D),\, (E\vert_{f^{-1}(D)})\otimes
N^*_{f^{-1}(D)})\, =\, 0\, ,
$$
where $N_{f^{-1}(D)}$ is the normal bundle of the
hypersurface $f^{-1}(D)$.
\end{proposition}

\begin{proof}
Consider the exact sequence of sheaves
$$
0\, \longrightarrow\, E\otimes {\mathcal O}_Z(-2f^{-1}(D))
\, \longrightarrow\,E\otimes {\mathcal O}_Z(-f^{-1}(D))
\, \longrightarrow\, (E\vert_{f^{-1}(D)})\otimes N^*_{f^{-1}(D)}
\, \longrightarrow\, 0
$$
over $Z$. Let
\begin{equation}\label{16o}
\longrightarrow\, H^1(Z,\, E\otimes {\mathcal O}_Z(-f^{-1}(D)))
\, \longrightarrow\, H^1(f^{-1}(D),\,
(E\vert_{f^{-1}(D)})\otimes N^*_{f^{-1}(D)})
\end{equation}
$$
\longrightarrow\, H^2(Z,\, E\otimes {\mathcal O}_Z(-2f^{-1}(D)))
$$
be the long exact sequence of cohomologies corresponding
to it.

Proposition \ref{vanish4}(1) says that
there is an integer $N'$ such that
$$
H^1(Z,\, E\otimes {\mathcal O}_Z(-f^{-1}(D)))\, =\, 0
$$
for all $r\, >\, N'$. Therefore, using Lemma
\ref{vanish6}, from Eqn. \eqref{16o} we conclude that
there is an integer $N$ such that
$$
H^1(f^{-1}(D),\, (E\vert_{f^{-1}(D)})\otimes N^*_{f^{-1}(D)})
\, =\, 0
$$
for all $r\, >\, N$. This completes the proof of the
proposition.
\end{proof}

\section{Injectivity of homomorphism of fundamental group
schemes}\label{sec.-inj.}

As before, let $X$ be a smooth projective variety defined
over an algebraically closed
field $k$ of characteristic $p$, with $p\, >\, 0$.
Let $L$ be an ample line bundle over $X$.

\begin{theorem}\label{thm2}
Assume that $\dim X \, \geq\, 3$. There is an integer
$d_0\, =\, d(X,L)$ with the following property.
Let $D\, \in\, \vert L^{\otimes d}\vert$ be a smooth
divisor, where $d\, >\, d_0$. Let $x_0$ be any $k$--rational
point of $D$. Then the
homomorphism between fundamental group schemes
$$
\pi(D,x_0) \,\longrightarrow\, \pi(X,x_0)
$$
induced by the inclusion map $D\, \hookrightarrow\, X$
is a closed immersion.
\end{theorem}

\begin{proof}
We will use the criterion given in Proposition
\ref{inj-surj}(2) for a homomorphism to be a
closed immersion.

Let
\begin{equation}\label{iD}
D\, \in\, \vert L^{\otimes d}\vert
\end{equation}
be a smooth divisor. Let $E$ be an essentially finite
vector bundle over $D$. Using Proposition
\ref{inj-surj}(2), to prove the theorem it suffices to
show the following:

If $d\, \geq\, d_0$, then
there is an essentially finite vector bundle
$V$ over $X$ whose restriction $V\vert_D$ to $D$ is
isomorphic to $E$.

Take any essentially finite vector bundle $E$ over $D$.
Let
\begin{equation}\label{i.e1}
f_D\, :\, \widetilde{D}\, \longrightarrow\, D
\end{equation}
be a connected \'etale Galois covering such that $f^*_D E$ is an
$F$--trivial vector bundle over $\widetilde{D}$ (see
Definition \ref{def1} for $F$--trivial vector bundles).

As $\dim X \, \geq\, 3$ with $D$ ample, from
Grothendieck's Lefschetz theory we know that the inclusion of $D$
in $X$ induces an isomorphism of \'etale fundamental groups
\cite[page 123, Th\'eor\`eme 3.10]{Gr},
\cite[page 177, Corollary 2.2]{Ha}. Therefore, there is
a unique (up to an isomorphism) connected \'etale Galois covering
\begin{equation}\label{i.e2}
f\, :\, Z\, \longrightarrow\, X
\end{equation}
and an isomorphism 
\begin{equation}\label{XYZ1}
\alpha\, :\, \widetilde{D}\,
\longrightarrow \, f^{-1}(D)
\end{equation}
such that $f_D\, =\,
f\circ\alpha$, where $f_D$ is the covering
morphism in Eqn. \eqref{i.e1}.

Let
\begin{equation}\label{i.e3}
F_{\widetilde{D}}\, :\, \widetilde{D}\, \longrightarrow\,
\widetilde{D}
\end{equation}
be the Frobenius morphism of the variety $\widetilde{D}$.
Since $f^*_D E$ is an $F$--trivial vector bundle, there
is a nonnegative integer $m$ such that the vector
bundle $(F^m_{\widetilde{D}})^* f^*_D E$
over $\widetilde{D}$ is trivializable.
As done in the proof of Theorem \ref{thm1}, we will employ
induction on $m$.

First assume that $m\, =\, 0$, in other words,
the vector bundle $f^*_D E$ itself is trivializable. As
$d\, \geq\, d'_0$, the inclusion map $D\,\hookrightarrow\,X$
gives an isomorphism of \'etale fundamental groups.
Since $f^*_D E$ is trivializable, this implies
that there is a vector bundle $V$ over $X$ whose restriction
$V\vert_D$ is isomorphic to $E$, and
furthermore, the vector bundle $f^* V$ is trivializable,
where $f$ is the morphism in Eqn. \eqref{i.e2}. In particular,
the criterion in Proposition \ref{inj-surj}(2) is valid if
$m\, =\, 0$.

We assume that $E$ extends to an essentially finite vector
bundle over $X$ provided $m\, \leq\, n_0-1$. We will show
that $E$ extends to an essentially finite vector
bundle over $X$ if $m\, =\, n_0$.

Let $E$ be an essentially finite vector bundle over $D$ with
$m\, =\, n_0$. Consider the vector bundle $F^*_D E$ over $D$,
where $F_D\, :\, D\, \longrightarrow\, D$ is the
Frobenius morphism. Since
\begin{equation}\label{i.is.}
f^*_D F^*_D E\, =\, F^*_{\widetilde{D}} f^*_D E
\end{equation}
over $\widetilde{D}$, where $f_D$ is the morphism in
Eqn. \eqref{i.e1} and $F_{\widetilde{D}}$ is defined
in Eqn. \eqref {i.e3}, we conclude that the vector bundle
$(F^{n_0-1}_{\widetilde{D}})^* f^*_D F^*_D E$ is
isomorphic to $(F^{n_0}_{\widetilde{D}})^* f^*_D E$.
Hence $(F^{n_0-1}_{\widetilde{D}})^* f^*_D F^*_D E$
is trivializable. Therefore, by the induction hypothesis,
the vector bundle $F^*_D E$ extends to $X$ as an essentially
finite vector bundle.

Let $V'$ be an essentially finite vector bundle over $X$
such that $V'\vert_D \, \cong\, F^*_D E$. Fix an
isomorphism of $V'\vert_D$ with $F^*_D E$. Let
\begin{equation}\label{iV1}
V_1\, :=\, f^* V'
\end{equation}
be the essentially finite vector bundle over $Z$,
where $f$ is the morphism in Eqn. \eqref{i.e2}.

Henceforth, we will assume that $d\, >\,
pL_{\rm max}(\Omega_X)/{\rm degree}(L)$.

We will identify $\widetilde D$ with $f^{-1}(D)$ using $\alpha$
in Eqn. \eqref{XYZ1}. In view of the above assumption
on $d$ and Eqn. \eqref{eqq1}, from Theorem \ref{thm1}
it follows that the homomorphism
of fundamental group schemes
\begin{equation}\label{i.sh}
\pi(\widetilde{D},x_0) \,\longrightarrow\, \pi(Z,x_0)
\end{equation}
induced by the inclusion map $\widetilde{D}\, \hookrightarrow\,
Z$ is surjective, where $x_0$ is any $k$--rational point
of $\widetilde D$. Note that using Eqn. \eqref{i.is.},
Eqn. \eqref{iV1} and the identity $\iota \circ F_{\widetilde{D}}\,
=\, F_Z\circ \iota$, where $\iota\, :\, \widetilde{D}\,
\longrightarrow\, Z$ is the inclusion map
and $F_Z$ is the Frobenius morphism of $Z$, we have
$$
((F^{n_0-1}_Z)^* V_1)\vert_{\widetilde{D}}\, =\,
(F^{n_0-1}_{\widetilde{D}})^*((f^*V')\vert_{\widetilde{D}})
\, =\, (F^{n_0-1}_{\widetilde{D}})^*f^*_D F^*_D E
\,=\, (F^{n_0}_{\widetilde{D}})^* f^*_D E\, .
$$
Hence $((F^{n_0-1}_Z)^* V_1)\vert_{\widetilde{D}}$
is trivializable as $(F^{n_0}_{\widetilde{D}})^* f^*_D E$ is
so (by assumption on $E$). Therefore, from the fact that
the homomorphism in Eqn. \eqref{i.sh} is surjective
it follows immediately that the vector bundle
$(F^{n_0-1}_Z)^* V_1$ over $Z$ is trivializable.

Using the isomorphism in Eqn. \eqref{i.is.}, the
Cartier connection on $F^*_{\widetilde{D}} f^*_D E$ induces
a connection on the vector bundle $f^*_D F^*_D E$. This
induced connection on $f^*_DF^*_DE$ 
will also be called the Cartier connection.
Let $\nabla^{\widetilde{D}}$ denote the Cartier connection on
$f^*_D F^*_D E$. This connection $\nabla^{\widetilde{D}}$
coincides with the pull
back, by $f_D$, of the Cartier connection on $F^*_D E$. Indeed,
this follows from the fact that the Cartier connection is
compatible with the pull back operation.

We will show that $\nabla^{\widetilde{D}}$ extends
to a connection on the vector bundle $V_1$ defined
in Eqn. \eqref{iV1}. For that we will first show that
$V_1$ admits a connection.

We recall from \cite{At} that a connection on $V_1$
is a splitting of the \textit{Atiyah exact sequence}
$$
0\, \longrightarrow\, \mathcal{E}nd(V_1)\, \longrightarrow\,
\text{At}(V_1) \, \longrightarrow\, TZ \, \longrightarrow\, 0\, .
$$
Let
\begin{equation}\label{i.obs.}
\theta(V_1)\, \in\, H^1(Z,\, \mathcal{E}nd(V_1)\otimes\Omega_Z)
\end{equation}
be the obstruction class for the existence of a splitting
of the Atiyah exact sequence. This obstruction class is
also called the \textit{Atiyah class} of $V_1$.

Let
\begin{equation}\label{i.ta.}
\tau\, :\, H^1(Z,\, \mathcal{E}nd(V_1)\otimes\Omega_Z)
\, \longrightarrow\,
H^1(\widetilde{D},\, (\mathcal{E}nd(V_1)\otimes
\Omega_Z)\vert_{\widetilde{D}})
\end{equation}
be the homomorphism induced by the inclusion map
$\widetilde{D}\, \stackrel{\iota}{\hookrightarrow}\, Z$.
For notational
convenience, the vector bundle $f^*_D F^*_D E \, =\,
V_1\vert_{\widetilde{D}}$ over $\widetilde{D}$
will be denoted by $E_1$. Let
\begin{equation}\label{i.ta.p}
\tau'\, :\, H^1(\widetilde{D},\, (\mathcal{E}nd(V_1)\otimes
\Omega_Z)\vert_{\widetilde{D}}) \, \longrightarrow\,
H^1(\widetilde{D},\, \mathcal{E}nd(E_1)\otimes
\Omega_{\widetilde{D}})
\end{equation}
be the homomorphism obtained from the natural projection
$$
\Omega_Z\vert_{\widetilde{D}}\, \longrightarrow\,
\Omega_{\widetilde{D}}\, .
$$

{}From general properties of the Atiyah class it
follows that
$$
c(E_1) \, :=\,\tau'\circ\tau (\theta(V_1))\, \in\,
H^1(\widetilde{D},\, \mathcal{E}nd(E_1)\otimes
\Omega_{\widetilde{D}})
$$
is the Atiyah class for $E_1$ where $\theta(V_1)$, $\tau$
and $\tau'$ are defined in Eqn. \eqref{i.obs.}, Eqn. \eqref{i.ta.}
and Eqn. \eqref{i.ta.p} respectively. Since $E_1$ admits
the Cartier connection $\nabla^{\widetilde{D}}$, 
the Atiyah class $c(E_1)$ vanishes.

Since $\tau'\circ\tau (\theta(V_1)) \, =\, c(E_1) \, =\, 0$,
we conclude that $\theta(V_1)\, =\, 0$ if both the
homomorphisms $\tau$ and $\tau'$ are injective.

Consider the long exact sequence of cohomologies
$$
H^1(Z,\, \mathcal{E}nd(V_1)\otimes\Omega_Z
(-\widetilde{D})) \, \longrightarrow\,
H^1(Z,\, \mathcal{E}nd(V_1)\otimes
\Omega_Z)\, \stackrel{\tau}{\longrightarrow}
\, H^1(\widetilde{D},\, \mathcal{E}nd(E_1)\otimes
(\Omega_Z\vert_{\widetilde{D}}))
$$
corresponding to the short exact sequence of sheaves
\begin{equation}\label{i.-lc1}
0\, \longrightarrow\, \mathcal{E}nd(V_1)\otimes\Omega_Z
(-\widetilde{D}) \, \longrightarrow\,\mathcal{E}nd(V_1)\otimes
\Omega_Z \, \longrightarrow\, \mathcal{E}nd(E_1)\otimes
(\Omega_Z\vert_{\widetilde{D}}) \, \longrightarrow\, 0\, .
\end{equation}
Using it together with Proposition \ref{vanish4}(2), which says that
there exists an integer $N(X,L)$ such that
\begin{equation}\label{i.-lc2}
H^1(Z,\, \mathcal{E}nd(V_1)\otimes\Omega_Z
(-\widetilde{D}))\, =\, 0
\end{equation}
if $r\, >\, N(X,L)$,
we conclude that the homomorphism $\tau$ in Eqn. \eqref{i.ta.}
is injective if $r\, >\, N(X,L)$ .

To prove that $\tau'$ is injective, consider the long
exact sequence of cohomologies
$$
H^1(\widetilde{D},\, \mathcal{E}nd(E_1)
\otimes N^*_{\widetilde{D}}) \, \longrightarrow\,
H^1(\widetilde{D},\, \mathcal{E}nd(E_1)\otimes
(\Omega_Z\vert_{\widetilde{D}}))
\, \stackrel{\tau'}{\longrightarrow}
\, H^1(\widetilde{D},\, \mathcal{E}nd(E_1)\otimes
\Omega_{\widetilde{D}})
$$
corresponding to the short exact sequence of sheaves
\begin{equation}\label{i.-lc3}
0\, \longrightarrow\,\mathcal{E}nd(E_1)\otimes N^*_{\widetilde{D}}
\, \longrightarrow\,\mathcal{E}nd(E_1)\otimes
(\Omega_Z\vert_{\widetilde{D}}) \, \longrightarrow\,
\mathcal{E}nd(E_1)\otimes \Omega_{\widetilde{D}}
\, \longrightarrow\, 0\, ,
\end{equation}
where $N_{\widetilde{D}}$ as before is the normal bundle
to the divisor $\widetilde{D}$. Using
it together with Proposition \ref{vi},
which says that there exists some integer $N'(X,L)$
such that
\begin{equation}\label{i.-lc4}
H^1(\widetilde{D},\, \mathcal{E}nd(E_1)
\otimes N^*_{\widetilde{D}})\, =\, 0
\end{equation}
if $r\, >\, N'(X,L)$,
we conclude that $\tau'$ defined in
Eqn. \eqref{i.ta.p} is also injective.

Henceforth, we will assume that $r\, >\, N(X,L), N'(X,L)$.

Since both $\tau$ and $\tau'$ are injective and $\tau'\circ
\tau (\theta(V_1)) \, =\, c(E_1) \, =\, 0$, the Atiyah class
$\theta(V_1)$ defined in Eqn. \eqref{i.obs.} vanishes.
In other words, the vector bundle $V_1$ admits a connection.

We will now show that $V_1$ admits a connection that restricts
to the Cartier connection $\nabla^{\widetilde{D}}$ on
$V_1\vert_{\widetilde{D}} \,=\, E_1$.

Let $\nabla^{V_1}$ be a connection on $V_1$. Let
$\nabla^1$ denote the connection on
$E_1\, =\, V_1\vert_{\widetilde{D}}$ induced by
this connection $\nabla^{V_1}$. On the other hand, the
vector bundle
$E_1$ has the Cartier connection $\nabla^{\widetilde{D}}$.
Therefore, we have
\begin{equation}\label{i.c.d}
\gamma(\nabla^1) \, :=\, \nabla^{\widetilde{D}}
- \nabla^1\, \in\, H^0(\widetilde{D},\, \mathcal{E}nd(E_1)
\otimes\Omega_{\widetilde{D}})\, .
\end{equation}

Consider the long exact sequence of cohomologies
$$
H^0(Z,\, \mathcal{E}nd(V_1)\otimes\Omega_Z)
\stackrel{\tau_0}{\longrightarrow}\,
H^0(\widetilde{D},\, \mathcal{E}nd(E_1)\otimes
(\Omega_Z\vert_{\widetilde{D}}))\, \longrightarrow\,
H^1(Z,\, \mathcal{E}nd(V_1)\otimes\Omega_Z(-\widetilde{D}))
$$
corresponding to the exact sequence in Eqn.
\eqref{i.-lc1}. In view of Eqn. \eqref{i.-lc2}, it follows from
this long exact sequence that the restriction homomorphism
\begin{equation}\label{141}
\tau_0\, :\, H^0(Z,\, \mathcal{E}nd(V_1)\otimes
\Omega_Z)\, \longrightarrow
\, H^0(\widetilde{D},\, \mathcal{E}nd(E_1)\otimes
(\Omega_Z\vert_{\widetilde{D}}))
\end{equation}
is surjective.

Similarly, using Eqn. \eqref{i.-lc4}
and the long exact sequence of cohomologies
\begin{equation}\label{142}
H^0(\widetilde{D},\, \mathcal{E}nd(E_1)\otimes
(\Omega_Z\vert_{\widetilde{D}}))\,
\stackrel{\tau'_0}{\longrightarrow} \, H^0(\widetilde{D},\,
\mathcal{E}nd(E_1)\otimes \Omega_{\widetilde{D}})\,\longrightarrow
\, H^1(\widetilde{D},\, \mathcal{E}nd(E_1)
\otimes N^*_{\widetilde{D}})
\end{equation}
corresponding to the exact sequence in Eqn. \eqref{i.-lc3} we
conclude that the above homomorphism $\tau'_0$
is surjective. Therefore, the composition homomorphism
\begin{equation}\label{XYZ2}
\tau'_0\circ\tau_0\, :\, H^0(Z,\,\mathcal{E}nd(V_1)\otimes
\Omega_Z)\,\longrightarrow\, H^0(\widetilde{D},\,
\mathcal{E}nd(E_1)\otimes \Omega_{\widetilde{D}})
\end{equation}
is surjective, where $\tau_0$ and $\tau'_0$ are defined in
Eqn. \eqref{141} and Eqn. \eqref{142} respectively.

Fix any
\begin{equation}\label{i.d.b5}
\beta\, \in\, H^0(Z,\,\mathcal{E}nd(V_1)\otimes
\Omega_Z)
\end{equation}
such that
\begin{equation}\label{150}
\tau'_0\circ\tau_0(\beta)\, =\,
\gamma(\nabla^1)\, ,
\end{equation}
where the section $\gamma(\nabla^1)$ is constructed
in Eqn. \eqref{i.c.d} and $\tau'_0\circ\tau_0$ is the
surjective homomorphism constructed in Eqn. \eqref{XYZ2}.

Let
\begin{equation}\label{143}
\nabla'\, :=\, \nabla^{V_1}+\beta
\end{equation}
be the connection on $V_1$, where $\beta$ is defined in Eqn.
\eqref{i.d.b5} and $\nabla^{V_1}$ is the earlier mentioned
connection on $V_1$ (recall that $\nabla^{V_1}$ restricts to
$\nabla^1$ on $\widetilde{D}$).
Since $\nabla^1$ is the restriction of
$\nabla^{V_1}$ to $\widetilde{D}$, from Eqn. \eqref{150}
it follows immediately that the restriction
of the connection $\nabla'$ (defined in Eqn. \eqref{143})
to $\widetilde{D}$ coincides with
the Cartier connection $\nabla^{\widetilde{D}}$ on $E_1$.

Let $\Gamma \, =\, \text{Gal}(f)$ be the Galois group for
the covering $f$ defined in Eqn. \eqref{i.e2}. The vector
bundle $V_1$ in Eqn. \eqref{iV1}, being a pull back by $f$,
is equipped with an action of $\Gamma$.

For any $g\in \text{Gal}(f)$ the connection $g^*\nabla'$ is another
connection on $V_1$.
Since the space of all connections on $V_1$ is an affine space for
the vector space $H^0(Z,\, \mathcal{E}nd(V_1)\otimes\Omega_Z)$, 
we conclude that $g^*\nabla'-\nabla'$ is
an element of $H^0(Z,\, \mathcal{E}nd(V_1)\otimes\Omega_Z)$. 

We showed that the restriction of the connection
$\nabla'$ to $\widetilde D$ coincides with the
Cartier connection $\nabla^{\widetilde D}$. On
the other hand, the Cartier connection
$\nabla^{\widetilde D}$ is invariant under the
action of
$$
\text{Gal}(f_D)\, =\, \text{Gal}(f)\, =\, \Gamma
$$
on $E_1$. Therefore, the element
$$
g^*\nabla' -\nabla'\, \in\, H^0(Z,\,
\mathcal{E}nd(V_1)\otimes\Omega_Z)
$$
is in the kernel of the homomorphism
$$
\tau'_0\circ\tau_0\, :\,
H^0(Z,\, \mathcal{E}nd(V_1)\otimes\Omega_Z)\, \longrightarrow \,
H^0(\widetilde{D},\,\mathcal{E}nd(E_1)\otimes \Omega_{\widetilde{D}})\, ,
$$
where $\tau'_0\circ\tau _0$ is the homomorphism in Eqn. \eqref{XYZ2}.
This homomorphism $\tau'_0\circ\tau_0$ coincides with the composition
of the homomorphism in Eqn. \eqref{145} with
the homomorphism $\beta$ in Eqn. \eqref{e4}
with the substitution $E'\, =\, V_1\, =\, V'$.
In the proof of Theorem \ref{thm1}
we saw that the homomorphism in Eqn.
\eqref{145} and the homomorphism $\beta$
in Eqn. \eqref{e4} are both injective.
Therefore, the homomorphism $\tau'_0\circ\tau_0$
is also injective.
We showed above that
$$
\tau'_0\circ\tau_0 (g^*\nabla'-\nabla')
\, =\, 0\, .
$$
Consequently, $g^*\nabla'
\, =\, \nabla'$.

Hence we have proved that the connection $\nabla'$ on $V_1$ is 
preserved by the action of $\Gamma$ on $V_1$. Therefore, it 
descends to a connection on the vector bundle
$V'$ over $X$ (see Eqn. \eqref{iV1}). Let $\nabla$ denote
the connection on $V'$ given by $\nabla'$. Therefore, we have
\begin{equation}\label{147}
\nabla'\, =\, f^*\nabla\, .
\end{equation}
Moreover, it follows that the
restriction to the divisor
$D$ of the connection $\nabla$ on $V'$ coincides
with the Cartier connection on $F^*_D E$ (recall that
$V'\vert_D \, =\, F^*_D E$).
Indeed, this is a consequence of the fact that
the restriction of the connection $\nabla'$ to
$\widetilde D$ coincides with the Cartier
connection $\nabla^{\widetilde D}$.

We will now show that the $p$--curvature of $\nabla$
vanishes (see \cite[page 190, (5.0.4)]{Ka} for the
definition of $p$--curvature of a connection).

Let
\begin{equation}\label{i.p.cu.}
\psi(\nabla)\, \in \, H^0(X,\, \mathcal{E}nd(V')
\otimes F^*_X\Omega_X)
\end{equation}
be the $p$--curvature of the connection $\nabla$ on $V'$
in Eqn. \eqref{147}. Let
$$
{\widetilde \alpha}\, :\,
H^0(X,\, \mathcal{E}nd(V')\otimes F^*_X\Omega_X)\, \longrightarrow
\, H^0(D,\, (\mathcal{E}nd(V')\otimes F_D^*\Omega_X)\vert_{D})
$$
be the restriction homomorphism. In view of our assumption that
$$
d\, >\, p\cdot L_{\rm max}(\Omega_X)/{\rm degree}(L)\, ,
$$
{}from the second part
of the Lemma \ref{vanish1} we know that
$$
H^0(X,\, \mathcal{E}nd(V')\otimes (F^*_X\Omega_X)(-D))\, =\, 0\, .
$$
Hence the above homomorphism ${\widetilde \alpha}$ is injective.
Let
\begin{equation}\label{00}
{\widetilde \beta}\, :\,
H^0(D,\, (\mathcal{E}nd(V')\otimes F^*_X\Omega_X)\vert_{D})
\, \longrightarrow\,
H^0(D,\, \mathcal{E}nd(V'\vert_D)\otimes F_D^*\Omega_{D})
\end{equation}
be the homomorphism induced by the natural projection
$\Omega_X\vert_{D}\, \longrightarrow\,
\Omega_{D}$. Since
$$
H^0(D,\, \mathcal{E}nd(V'\vert_D)\otimes F_D^*N^*_{D})
\, =\, 0
$$
(see the second part of Lemma \ref{vanish2}), from the left exact
sequence of global sections for the short exact sequence
of sheaves
$$
0\, \longrightarrow\,\mathcal{E}nd(V'\vert_D)\otimes F_D^*N^*_{D}
\, \longrightarrow\,(\mathcal{E}nd(V')\otimes F^*_X\Omega_X)\vert_{D}
\,\longrightarrow\, \mathcal{E}nd(V'\vert_D)\otimes F_D^*\Omega_{D}
\, \longrightarrow\, 0
$$
it follows that the homomorphism ${\widetilde \beta}$
defined in Eqn. \eqref{00} is injective.

Since the $p$--curvature of any Cartier connection vanishes, and the
restriction to $D$ of the connection $\nabla$ coincides
with the Cartier connection on $F^*_D E$, we have
$$
{\widetilde \beta}\circ{\widetilde \alpha}(\psi(\nabla))\, =\, 0\, ,
$$
where $\psi(\nabla)$ is the $p$--curvature in Eqn. \eqref{i.p.cu.}.
As both ${\widetilde \alpha}$ and ${\widetilde \beta}$ are
injective homomorphisms, this implies that $\psi(\nabla)\, =\, 0$.

Since the $p$--curvature of $\nabla$ vanishes, there
is a vector bundle $W$ over $X$ such that
$F^*_X W$ equipped with the Cartier connection is identified
with the vector bundle $V'$ equipped with the connection $\nabla$ 
\cite[page 190, Theorem 5.1.1]{Ka}.

Since the restriction of $\nabla$ to $D$ coincides
with the Cartier connection on $F^*_D E$, the isomorphism
of $(F^*_X W)\vert_D$ with $V'\vert_D$ is the pull back, by
the Frobenius morphism $F_D$, of an isomorphism of $E$ with
$W\vert_D$. Since $F^*_XW \, =\, V'$ is an essentially finite vector
bundle, it follows that the vector bundle $W$ is also essentially finite
(this was noted in the proof of Theorem \ref{thm1}).

Therefore, using induction, any essentially finite
vector bundle over $D$ extends to $X$ as an essentially finite
vector bundle. Using Proposition \ref{inj-surj}(2)
this completes the proof of the theorem.
\end{proof}

Using the proof of Theorem \ref{thm2}
we have the following corollary:

\begin{corollary}\label{pn}
Let $D$ be a smooth hypersurface in
the projective space ${\mathbb P}^n_k$,
with $n\, \geq\, 3$. Then
the fundamental group scheme $\pi(D,x_0)$
is trivial.
\end{corollary}

\begin{proof}
We will follow the steps of the proof of
Theorem \ref{thm2}.

Take any smooth hypersurface $D$ of
${\mathbb P}^n_k$.
Since the \'etale fundamental group of $D$ is
trivial, we have $Z\, =\, {\mathbb P}^n_k$ and
$f$ to be the identity map of ${\mathbb P}^n_k$.
Next we note that the fundamental group scheme
of ${\mathbb P}^n_k$ is trivial 
(see the corollary in \cite[page 93]{No2} following
the proof of Proposition 9). 
In other words, any essentially finite
vector bundle over ${\mathbb P}^n_k$ is a
trivializable vector bundle. Next we note
that the cotangent bundle
$\Omega_{{\mathbb P}^n_k}$ is strongly semistable.
This can be proved using the facts that the
tangent bundle $T_{{\mathbb P}^n_k}$ is globally
generated and $\text{Aut}({\mathbb P}^n_k)$
acts irreducibly on $H^0({\mathbb P}^n_k,\,T_{{\mathbb P}^n_k})$.

In particular,
$$
L_{\rm max}(\Omega_{{\mathbb P}^n_k})\, =\,
- \frac{n+1}{n}\, .
$$

Using the above observations in the proof
of Theorem \ref{thm2} we conclude
that the proof of the corollary will be
complete once we establish the following two
assertions:
\begin{equation}\label{fich5p}
H^1({\mathbb P}^n_k,\, \Omega_{{\mathbb P}^n_k}
(-D))\, =\, 0
\end{equation}
and
\begin{equation}\label{sech5p}
H^1(D,\, N^*_D) \, =\, 0\, .
\end{equation}
Note that
these two correspond to Eqn. \eqref{i.-lc2}
and Eqn. \eqref{i.-lc4} respectively.

Since $n\, \geq\, 3$, we have
$$
H^1({\mathbb P}^n_k,\, {\mathcal O}_{{\mathbb
P}^n_k}(\delta)) \, =\, 0\, =\, H^2({\mathbb
P}^n_k,\, {\mathcal O}_{{\mathbb P}^n_k}(\delta))
$$
for all $\delta$. Therefore, Eqn. \eqref{fich5p}
follows from the long exact sequence of cohomologies
for the short exact sequence of vector bundles
$$
0\, \longrightarrow\, \Omega_{{\mathbb P}^n_k}(-D)
\, \longrightarrow\, {\mathcal O}_{{\mathbb
P}^n_k}(-1)(-D)^{\oplus (n+1)}\, \longrightarrow\,
{\mathcal O}_{{\mathbb P}^n_k}(-D)\,
\longrightarrow\, 0
$$
obtained from the Euler sequence.

Similarly, Eqn. \eqref{sech5p} follows from the
long exact sequence of cohomologies for the short
exact sequence
$$
0\, \longrightarrow\, {\mathcal O}_{{\mathbb
P}^n_k}(-2D)\, \longrightarrow\,
{\mathcal O}_{{\mathbb P}^n_k}(-D)\,\longrightarrow
\, N^*_D\, \longrightarrow\, 0\, .
$$
This completes the proof of the corollary.
\end{proof}

\begin{remark}\label{inj.eff.}
{\rm The integer $d_0$ in Theorem \ref{thm2} can be taken
to be any integer satisfying the following conditions.}
\begin{enumerate}
\item{} {\rm $d_0\, \geq\, p{\cdot}r_1$, where $r_1$ is
defined in Lemma \ref{etalevanish}.}

\item{} {\rm $d_0\, \geq \, N(B_1)$, where $B_1$ is defined in
Eqn. \eqref{eq55}, and $N(V)$ of a vector bundle $V$
is defined in the proof of Lemma \ref{vanish3}.}

\item{} {\rm $d_0\, \geq \, N({\mathcal G}_1, \cdots ,
{\mathcal G}_{(p-1)\dim (X)})$, where ${\mathcal G}_i$ are
defined in Proposition \ref{canfil}(3) and $N(G_1, \cdots ,
G_{(p-1)\dim (X)})$
is defined in the proof of Lemma \ref{vanish3-1}.}
\end{enumerate}
\end{remark}


\end{document}